\documentclass[11pt,a4paper]{amsart}
\usepackage[T1]{fontenc}
\usepackage{lmodern}
\usepackage{amsmath,amssymb,amsthm,mathtools}
\usepackage[margin=28mm]{geometry}
\usepackage{microtype}
\usepackage{booktabs,array}
\usepackage[colorlinks=true,linkcolor=blue!50!black,citecolor=blue!50!black,urlcolor=blue!50!black]{hyperref}
\usepackage{xcolor}
\hypersetup{pdftitle={Evasive subspaces, minimal rank-metric codes and subspace designs via auxiliary fields},pdfauthor={Alessandro Giannoni}}
\newcommand{\dimq}{\dim_{\mathbb F_q}}
\newcommand{\gauss}[3]{\genfrac{[}{]}{0pt}{}{#1}{#2}_{#3}}
\DeclareMathOperator{\rk}{rk}
\DeclareMathOperator{\wt}{wt}

\DeclareMathOperator{\Supp}{Supp}
\DeclareMathOperator{\Tr}{Tr}

\DeclareMathOperator{\PG}{PG}
\newtheorem{theorem}{Theorem}[section]
\newtheorem{lemma}[theorem]{Lemma}
\newtheorem{proposition}[theorem]{Proposition}
\newtheorem{corollary}[theorem]{Corollary}
\theoremstyle{definition}

\newtheorem{example}[theorem]{Example}
\theoremstyle{remark}
\newtheorem{remark}[theorem]{Remark}
\numberwithin{equation}{section}
\allowdisplaybreaks[2]
\title[Evasive subspaces, minimal codes and designs]{Evasive subspaces, minimal rank-metric codes\\and subspace designs via auxiliary fields}
\author{Alessandro Giannoni}
\address{Universit\`a degli Studi di Napoli Federico II, Naples, Italy}
\email{alessandro.giannoni@unina.it}
\begin{document}
\begin{abstract}
We extend the auxiliary-field method for scattered subspaces to evasive
subspaces and subspace designs. For $1\le h\le k<hn$ and coprime $h,k+1$
with $n(k-h+1)\ge k+1$, we construct spanning $(h,k)$-evasive
$\mathbb F_q$-subspaces of $\mathbb F_{q^n}^{\,r}$ of dimension
$rn(k-h+1)/(k+1)$ in every sufficiently large admissible ambient
dimension, with a threshold independent of $q$. When $k<n$, this
dimension is maximum. The proof combines an intersection lemma,
a rank-support density inequality, scalar orbit counting, and a
profile preserved by direct sums. The coding consequences include,
for every fixed odd extension degree $m\ge3$ and every sufficiently
large $K$, minimal rank-metric codes of optimal length $2K$ with
first two generalized weights $(m+1)/2$ and $m+1$. For these codes the
stable profile is equivalent to minimality with the largest possible
minimum distance; for even $m$ the analogous class consists of minimal
MRD codes. More generally, for $1\le \tau\le m-2$ and
$\gcd(m,\tau+1)=1$, we obtain nondegenerate $\tau$-minimal codes of
length $(\tau+1)K$ for every sufficiently large $K$, uniformly in $q$.
We also construct maximum $h$-designs whose number of members grows
exponentially with the ambient dimension at fixed $q,n,h$, and the
associated maximum sum-rank distance (MSRD) codes with a growing
number of equal long blocks.
In extension degree five, the ordinary and $2$-minimal constructions
hold for every $q$ and every $K\ge32$ and $K\ge59$, respectively.
\end{abstract}
\maketitle
\section{Introduction}\label{sec:introduction}

Restricted intersections of $\mathbb F_q$-subspaces with subspaces
over $\mathbb F_{q^n}$ connect finite geometry with rank-metric
coding theory. Scattered spaces arose in the study of spreads and
linear blocking sets \cite{BBL,BL}; linear sets provide the
projective setting for these questions \cite{Polverino}.
Their relation with maximum rank-distance (MRD) codes, originating
with Delsarte and Gabidulin \cite{Delsarte,Gabidulin}, has led to
both extremal existence results and new constructions
\cite{BGMP,CMPZ17}. Csajb\'ok, Marino, Polverino and Zullo
\cite{CMPZ21} introduced $h$-scattered subspaces, with further
geometric and coding developments in \cite{ZZ,BGGM,BGM,BGGMN}.
The language of $q$-systems and generalized rank weights makes
these connections precise \cite{SVV,Rand,MNT}.

Bartoli, Csajb\'ok, Marino and Trombetti \cite{BCMT} developed the
more general notion of $(h,k)$-evasivity and its dimension bounds.
It belongs to a broader intersection problem: given a family of
subspaces, find a large space meeting every member in bounded
dimension. Gruica, Ravagnani, Sheekey and Zullo \cite{GRSZ}
formulate this as a generalised critical problem and study several
families by counting and density arguments. Related counting
questions for rank-metric codes also appear through rank-metric
lattices \cite{CR,CRZ}. Subspace evasive sets, without the linearity
requirement imposed here, were studied by Dvir and Lovett
\cite{DL} in connection with list decoding.

Giannoni and Neri \cite{GN} recently gave a probabilistic
construction of maximum $h$-scattered subspaces. Their
auxiliary-field method realizes the desired space as the image
of an $\mathbb F_{q^s}$-linear domain under a random matrix over
$\mathbb F_{q^{ns}}$. When $\gcd(n,s)=1$, scalar orbits reduce
the number of bad relations in the first-moment estimate.
They also use a finite \'etale tensor algebra to remove this
restriction on the auxiliary degree.

We extend the coprime-field method to general evasivity and to
families satisfying a joint intersection bound. The new obstruction
is that a bad relation is now an $\mathbb F_{q^n}$-linear space,
whose span over $\mathbb F_{q^{ns}}$ may lose dimension.
An intersection lemma and a rank-support density inequality control
this loss, while relation spaces whose span does not drop are counted
by scalar orbits, for the main constraint as well as for the additional
constraints of a profile. A profile preserved by direct sums then gives every
sufficiently large admissible ambient dimension. The scattered
case $h=k\le n-2$ is already covered by \cite{GN}, with a sharper
threshold, and $h=k=n-1$ is elementary. Our contributions concern
$2\le h<k$, the simultaneous
profiles, their minimal-code consequences, and designs with
many members. In particular, the one-member specialization of
our design criterion recovers the coprime-field criterion of
\cite{GN}.

\subsection{Main results}
A spanning $\mathbb F_q$-subspace
$U\le\mathbb F_{q^n}^{\,r}$ is $(h,k)$-evasive if every
$h$-dimensional $\mathbb F_{q^n}$-subspace meets it in dimension
at most $k$ over $\mathbb F_q$.

\begin{theorem}\label{thm:intro-evasive}
Fix integers $n\ge2$ and $1\le h\le k<hn$ such that
\[
 \gcd(h,k+1)=1,\qquad n(k-h+1)\ge k+1.
\]
There is an integer $r_0(n,h,k)$ such that, for every prime power
$q$ and every $r\ge r_0(n,h,k)$ with $(k+1)\mid rn$, there is
a spanning $(h,k)$-evasive subspace
\[
 U\le\mathbb F_{q^n}^{\,r},\qquad
 \dimq U=\frac{rn(k-h+1)}{k+1}.
\]
It can be chosen to satisfy
\[
 \dimq(U\cap W)\le
 \left\lfloor\frac{(k+1)d}{h}\right\rfloor
 \qquad\left(1\le d=\dim_{\mathbb F_{q^n}}W<h\right).
\]
When $k<n$, this dimension is maximum.
\end{theorem}

Theorem~\ref{thm:all-r} gives a computable threshold independent
of $q$. Maximality for $k<n$ follows from \cite{BCMT}.
For $n\le k<hn$ the result supplies a specific attainable
dimension; a general bound in this range, which follows from the
cardinality bound \cite[Theorem~4.3]{BCMT}, is
$\dimq U\le n(r-h)+k$ \cite[Corollary~4.6(2)]{MNT}.
When $k\ge hn$, evasivity is automatic, so this range is omitted
from the introductory statement.
Our result addresses attainment for the family of all
$\mathbb F_{q^n}$-linear $h$-spaces, rather than the arbitrary
families in the generalised critical problem of \cite{GRSZ}.
For $h=1$, that is, when every point of $\PG(r-1,q^n)$ must have
weight at most $k$, maximum evasive subspaces are already available;
duals are taken with respect to $(x,y)\mapsto\Tr_{q^n/q}(x\cdot y)$.
For every $k<n$, the dual of the maximum $k$-scattered subspace
$\{(x,x^q,\ldots,x^{q^k}):x\in\mathbb F_{q^n}\}$ of $\mathbb F_{q^n}^{\,k+1}$
is a spanning maximum $(1,k)$-evasive subspace, and direct sums give
every $r$ with $(k+1)\mid r$. For $k=n-2$ these are all admissible
ranks, and for $k=n-1$ the dual of $\mathbb F_q^r$ works for every $r$.
For $k\le n-3$, the remaining sufficiently large admissible ranks follow
from \cite{GN} and the duality \cite[Theorem~2.15]{GRSZ} between
maximum $k$-scattered and maximum $(1,k)$-evasive subspaces;
\cite[Corollary~4.1]{GRSZ} lists explicit cases. For arbitrary spreads
and $q\ge64^{1/(k+1)}$, \cite[Proposition~4.15]{GRSZ} gives evasive
subspaces of every dimension up to $rnk/(k+1)-k(n-k-1)/(k+1)$. Thus
Theorem~\ref{thm:intro-evasive} is new for $2\le h<k$, together with
the stable profile.
Theorem~\ref{thm:simultaneous} allows additional constraints
with strictly smaller density ratio.

The rank-weight dictionary of \cite{MNT} yields applications to
minimal codes. These were introduced and studied through linear
cutting blocking sets by Alfarano, Borello, Neri and Ravagnani
\cite{ABNR}. Following the literature on minimal codes, we denote the
extension degree by $m$ in this context, and write
$\varpi_{q^m/q}(K)$ for the minimum length
of an $\mathbb F_{q^m}$-linear minimal rank-metric code of
dimension $K$.

\begin{theorem}\label{thm:intro-minimal}
For each fixed odd $m\ge3$, there is $K_0(m)$ such that
\[
 \varpi_{q^m/q}(K)=2K
\]
for every prime power $q$ and every $K\ge K_0(m)$.
The codes attaining this length can be chosen with
\[
 d_1=\frac{m+1}{2},\qquad d_2=m+1.
\]
\end{theorem}

The length and weight assertions are
Theorem~\ref{thm:minimal-optimal} and
Corollary~\ref{cor:minimal-exact-weights}.
For $m=5$, Corollary~\ref{cor:minimal-degree-five} gives
$\varpi_{q^5/q}(K)=2K$, with $d_1=3$, for every $K\ge32$ and ten
smaller values of $K$, and also $d_2=6$ for every $K\ge34$;
Remark~\ref{rem:thresholds-odd} describes the
admissible seeds, and hence explicit thresholds, for every odd $m$.
The minimum-length problem is explicitly raised in
\cite[Remark 6.12]{ABNR} and discussed in
\cite[Section 4.3]{Scotti}.
The bounds relevant to our result are collected below;
all rows hold for every prime power $q$.

\begin{center}
\small
\begin{tabular}{@{}lll@{}}
\toprule
Bound for $\varpi_{q^m/q}(K)$ & Range & Source\\
\midrule
$\ge K+m-1$ & $m,K\ge2$ & \cite[Corollary 5.10]{ABNR}\\
$\ge 2K$ & $m\ge3,\ K\ge2$ & Proposition~\ref{prop:minimal-lower}; \cite[Theorem~5.3(1)]{XKH} for $m=3$\\
$\le 2K+m-2$ & $m,K\ge2$ & \cite[Theorem 6.11]{ABNR}\\
$\le 2K+m-3$ & $m\ge3,\ K\ge2$ & \cite[Corollary~5.4(1)]{XKH}\\
$=2K$ & odd $m\ge3,\ K\ge K_0(m)$ &
Theorem~\ref{thm:minimal-optimal}\\
\bottomrule
\end{tabular}
\end{center}

Xu, Kan and Han \cite[Corollary~5.4(2)]{XKH} already determine the
value $2K$ for $m=3$ and every $K\ge2$. Other short constructions include
the $[8,4]_{q^4/q}$ cutting codes of Bartoli, Marino and Neri
\cite{BMN} for $q$ an odd power of $2$, and the
$[m+2,3]_{q^m/q}$ codes of Lia, Longobardi, Marino and Trombetti
\cite{LLMT} for odd $m\ge5$ and infinitely many $q$.
Thus our new length conclusion is the exact eventual value
for every fixed odd $m\ge5$, uniformly in $q$. In this range it
also settles the existence question for the intermediate lengths
listed among the open problems in \cite[Section~5.2]{Scotti}:
minimal codes of dimension $K$ exist exactly for the lengths
$N\ge2K$, and nondegenerate ones for $2K\le N\le Km$ (see the
discussion after Theorem~\ref{thm:minimal-optimal}).

The profile used in the construction has a direct coding meaning.
For $[2K,K]_{q^m/q}$ codes with odd $m\ge5$, it is equivalent to
minimality together with the largest possible minimum distance
$(m+1)/2$, and in a minimal direct sum of $[2K_i,K_i]_{q^m/q}$ codes
every block has this property (Proposition~\ref{prop:profile-coding}).
For even $m$, the corresponding blocks are minimal MRD codes, and the
MRD condition has the same critical density ratio as minimality
(Proposition~\ref{prop:even-degree}). Combined with the
$[8,4]_{q^4/q}$ codes of \cite{BMN}, this gives
$\varpi_{q^4/q}(4w)=8w$ for every $w\ge1$ when $q$ is an odd power of
$2$ (Corollary~\ref{cor:quartic}).

Higher minimality compares the supports of $\tau$-dimensional
subcodes. Xu, Kan and Han \cite{XKH}, who call it $r$-minimality,
characterize it by $d_{\tau+1}\ge m\tau+1$ \cite[Theorem~4.4]{XKH}.
With $\varpi_{q^m/q}(K,\tau)$ denoting the minimum length of a
$\tau$-minimal code of dimension $K$, they give the general existence
bound
$\varpi_{q^m/q}(K,\tau)\le(\tau+1)K+\tau(m-\tau-2)$ for
$1\le \tau\le m-2$ and $K\ge \tau+1$ \cite[Theorem~5.4]{XKH}.

\begin{theorem}\label{thm:intro-higher-minimal}
Fix $m\ge3$ and $1\le \tau\le m-2$ with $\gcd(m,\tau+1)=1$.
For every sufficiently large $K$, with a threshold independent of
$q$, there is a nondegenerate $\tau$-minimal code with parameters
\[
 [(\tau+1)K,K]_{q^m/q}
\]
for every prime power $q$.
\end{theorem}

Theorem~\ref{thm:higher-minimal} removes the additive term
$\tau(m-\tau-2)$ under these hypotheses; this is a strict improvement
when $\tau\le m-3$. At $\tau=m-2$ the length is already attained
by \cite[Theorem~5.4]{XKH}. Corollary~\ref{cor:higher-minimal-five} gives
nondegenerate $2$-minimal $[3K,K]_{q^5/q}$ codes for every
$q$ and every $K\ge59$.

Subspace designs arose in list decoding
\cite{GX,GK,GWX,GXY}. Their finite-geometric form and its connection
with the sum-rank metric were developed in
\cite{SZ,SS}, alongside the geometric description of sum-rank
codes in \cite{NSZ}.
Here an $h$-design is a jointly spanning family
$(U_1,\ldots,U_t)$ in $\mathbb F_{q^n}^{\,r}$ satisfying
\[
 \sum_{i=1}^t\dimq(U_i\cap H)\le h
 \qquad\left(\dim_{\mathbb F_{q^n}}H=h\right).
\]
It is maximum when each member has dimension $rn/(h+1)$.

\begin{theorem}\label{thm:intro-design}
Fix $1\le h<n$. For every fixed positive integer $t$, a maximum
$h$-design with $t$ members exists in $\mathbb F_{q^n}^{\,r}$
for every prime power $q$ and every sufficiently large $r$ with
$(h+1)\mid rn$; the threshold is independent of $q$.
For fixed $q,n,h$, put $\nu=(h+1)/\gcd(n,h+1)$. For every
sufficiently large admissible $r$, the number of members can
be chosen to satisfy
\[
 t\ge
 \begin{cases}
 q^{\,r/(2\nu(h+1))-O_{q,n,h}(1)},& n\ \text{odd},\\
 q^{\,r/(3\nu(h+1))-O_{q,n,h}(1)},& n\ \text{even}.
 \end{cases}
\]
If $\gcd(r/\nu,n)=1$, the denominator is $\nu(h+1)$.
\end{theorem}

Theorems~\ref{thm:design-all-r} and \ref{thm:design-growth}
prove these assertions; the latter also gives the denominator
$2\nu(h+1)$ when $n$ is a power of $2$.
For $h=1$, the members define pairwise disjoint maximum
scattered linear sets.
For $t\le q-1$, maximum $h$-designs were known when $(h+1)\mid r$
\cite{SZ} (see also \cite[Theorem~6.7]{SS}) and, for $h=1$, whenever
$rn$ is even \cite[Theorem~4.11]{SS}. For $h=1$, families with
$2(q-1)$ members occur for $n=4$ in \cite[Corollary~4.15]{SS}, and for
$n=2$ and odd $r$ partitions into Baer subgeometries give
$(q^r+1)/(q+1)$ members \cite[Theorem~4.7]{SS}.
Theorem~\ref{thm:intro-design} gives any prescribed number of members
for every $1\le h<n$, in all sufficiently large admissible dimensions,
and exponential growth at fixed $q,n,h$; for $h=1$ it answers the
first question of \cite[Problem~4.17]{SS} in that range. For $n=2$,
$h=1$ and odd $r$, the Baer partitions remain much larger.
Corollary~\ref{cor:msrd} uses the correspondence of
\cite[Theorem~3.16]{SS} to obtain maximum sum-rank distance (MSRD)
codes with equal long blocks.
In particular, it gives such codes with block length $rn/(h+1)$ for
every $1\le h<n$ and every sufficiently large admissible $r$, as asked
in \cite[Problem~6.11]{SS}, although not by an explicit construction.
Section~\ref{sec:designs} places these codes relative to the
constructions in \cite{MP18,MP22,MP24}, the matrix-metric framework
of \cite{BGR}, and the recent length bounds of \cite{AZ}.

The constructions are existence arguments with explicit finite
criteria. Optimal length is proved for ordinary minimality in
the stated odd-degree range; for higher minimality we prove an
upper bound. The exact first two weights concern the selected
codes, and neither the numerical thresholds nor the design
exponents are claimed to be optimal.

\subsection{Structure of the argument}
After the preliminaries, we prove the intersection and density
lemmas. The finite existence criterion separates relation spaces
by the dimension of their span over $\mathbb F_{q^{ns}}$.
Scalar orbits control the full-span case, for the main constraint
as well as for the additional constraints of a profile, while strict density
inequalities control the smaller spans. Simultaneous constraints
and stable profiles then give all large admissible ranks.
The final sections develop generalized weights and minimality,
where the stable profile is identified with a coding property,
followed by maximum designs and MSRD codes.

\section{Geometric and coding preliminaries}\label{sec:preliminaries}

Throughout, $q$ is a prime power and $n$ is a positive integer. Dimensions over $\mathbb F_q$ are written $\dimq$. For a subset $X$ of an $\mathbb F_{q^n}$-space, $\langle X\rangle_{\mathbb F_{q^n}}$ denotes its $\mathbb F_{q^n}$-linear span. All subspaces contain zero.

\subsection{Evasivity and extremal dimension}
Let $1\le h<r$ and $h\le k$. An $\mathbb F_q$-subspace $U\le \mathbb F_{q^n}^r$ is \emph{$(h,k)$-evasive} if
\begin{equation}\label{eq:evasive-definition}
 \dimq(U\cap H)\le k
 \qquad\text{for every }H\le_{\mathbb F_{q^n}} \mathbb F_{q^n}^r\text{ with }\dim_{\mathbb F_{q^n}} H=h.
\end{equation}
We shall require $\langle U\rangle_{\mathbb F_{q^n}}=\mathbb F_{q^n}^r$ when calling $U$ a \emph{$q$-system}, following the geometric description of rank-metric codes in \cite{SVV,Rand,MNT}. This spanning convention is stronger than the convention in \cite{BCMT}, which only requires the span to have dimension at least $h$.

When $k<n$, we use the following dimension bound of
Bartoli, Csajb\'ok, Marino and Trombetti \cite[Corollary~4.9]{BCMT}.
\begin{theorem}\label{thm:evasive-bound}
If $1\le h\le k<n$, $r\ge h+1$, and $U\le \mathbb F_{q^n}^r$ is a spanning $(h,k)$-evasive $\mathbb F_q$-subspace, then
\[
 \dimq U\le \frac{rn(k-h+1)}{k+1}.
\]
\end{theorem}
A spanning evasive subspace attaining this rational bound is called \emph{maximum} in this paper. The distinction from attaining the integer part of the bound matters when the displayed quantity is not integral. The special case $k=h$ is the notion of an \emph{$h$-scattered} subspace \cite{CMPZ21}. When $h<n$, the corresponding dimension bound is $rn/(h+1)$.
For arbitrary $k$, the dimension inequality
$\dimq(U\cap H)\ge\dimq U+nh-nr$ also gives
$\dimq U\le n(r-h)+k$; see \cite[Theorem~4.3]{BCMT} and
\cite[Corollary~4.6(2)]{MNT}.

For $a=k+1$ and $j=a-h$, the coprimality assumption used below can be written in either of the forms
\[
 \gcd(h,a)=1\qquad\text{or}\qquad\gcd(j,a)=1.
\]
Under this assumption, integrality of $rnj/a$ is equivalent to $a\mid rn$.

\subsection{Rank supports and generalized weights}
We use the support description of generalized rank weights
\cite{KMU,Rav,MNT}. Fix an $\mathbb F_q$-basis
$e_1,\ldots,e_n$ of $\mathbb F_{q^n}$.
For $v\in \mathbb F_{q^n}^N$, write
\[
 v=\sum_{b=1}^n e_b v_b,\qquad v_b\in\mathbb F_q^N.
\]
The \emph{rank support} of $v$ is the $\mathbb F_q$-span of the rows $v_b$, and its dimension is $\rk_q(v)$. For an $\mathbb F_{q^n}$-linear subspace $D\le \mathbb F_{q^n}^N$, put
\[
 \Supp_q(D)=\sum_{v\in D}\Supp_q(v),\qquad
 \wt_q(D)=\dimq\Supp_q(D).
\]
These definitions are independent of the chosen basis. The generalized rank weights of an $[N,K]_{q^n/q}$ code $D$, meaning an $\mathbb F_{q^n}$-linear subspace of $\mathbb F_{q^n}^N$ of dimension $K$, are
\[
 d_j(D)=\min\{\wt_q(W):W\le_{\mathbb F_{q^n}} D,\ \dim_{\mathbb F_{q^n}} W=j\},
 \qquad 1\le j\le K.
\]
Equivalently,
\begin{equation}\label{eq:support-characterization}
 d_j(D)=\min\{\dimq B:B\le_{\mathbb F_q}\mathbb F_q^N,\
                   \dim_{\mathbb F_{q^n}}(D\cap \langle B\rangle_{\mathbb F_{q^n}})\ge j\}.
\end{equation}
Here $\langle B\rangle_{\mathbb F_{q^n}}$ is the $\mathbb F_{q^n}$-span of $B$. Indeed, $W\subseteq \langle B\rangle_{\mathbb F_{q^n}}$ is equivalent to $\Supp_q(W)\subseteq B$. We will use the generalized Singleton inequality
\begin{equation}\label{eq:generalized-singleton}
 d_j(D)\le N-K+j.
\end{equation}
One direct proof is to intersect $D$ with a coordinate subspace of dimension $N-K+j$, which leaves dimension at least $j$.

The code $D$ is \emph{nondegenerate} if $\Supp_q(D)=\mathbb F_q^N$. Equivalently, the columns of any generator matrix are $\mathbb F_q$-independent. With respect to the ordinary $\mathbb F_{q^n}$-bilinear dot product,
\[
 D^\perp=\{v\in \mathbb F_{q^n}^N:\ v\cdot w=0\text{ for all }w\in D\}.
\]
A code is nondegenerate if and only if its dual contains no nonzero vector of rank one. This follows because a rank-one vector is an $\mathbb F_{q^n}$-multiple of a nonzero vector of $\mathbb F_q^N$.

Choose an $\mathbb F_q$-basis of a spanning $q$-system $U\le \mathbb F_{q^n}^r$ as the columns of an $r\times N$ matrix $G$, where $N=\dimq U$. The row space $C_U$ is a nondegenerate $[N,r]_{q^n/q}$ code. Conversely, this construction recovers a $q$-system from a nondegenerate code, up to the usual choices of bases.

We record \cite[Theorem 3.3, $(1)\Leftrightarrow(3)$]{MNT}
in the notation needed below; the short proof fixes our conventions.
\begin{proposition}\label{prop:dictionary}
Let $U\le \mathbb F_{q^n}^r$ be a spanning $q$-system of $\mathbb F_q$-dimension $N$, let $G$ be as above, and put $D=\ker_{\mathbb F_{q^n}} G=C_U^\perp$. Suppose
\[
 1\le h<r,\qquad a=k+1\le N,\qquad
 1\le j=a-h\le N-r.
\]
Then $U$ is $(h,k)$-evasive if and only if $d_j(D)\ge k+2$.
\end{proposition}
\begin{proof}
The map $G:\mathbb F_q^N\to U$ is an $\mathbb F_q$-isomorphism. If an $h$-dimensional $\mathbb F_{q^n}$-space contains an $a$-dimensional $\mathbb F_q$-subspace of $U$, let $B\le\mathbb F_q^N$ be its inverse image. Then $\dim_{\mathbb F_{q^n}} \langle B\rangle_{\mathbb F_{q^n}}=a$ and $\dim_{\mathbb F_{q^n}} G(\langle B\rangle_{\mathbb F_{q^n}})\le h$, so $\dim_{\mathbb F_{q^n}}(D\cap \langle B\rangle_{\mathbb F_{q^n}})\ge a-h=j$. Formula \eqref{eq:support-characterization} gives $d_j(D)\le a$.

Conversely, if $d_j(D)\le a$, choose a support space $B$ as in \eqref{eq:support-characterization} and enlarge it to dimension $a$. The image $G(\langle B\rangle_{\mathbb F_{q^n}})$ has $\mathbb F_{q^n}$-dimension at most $a-j=h$, whereas $G(B)$ has $\mathbb F_q$-dimension $a$. Extending $G(\langle B\rangle_{\mathbb F_{q^n}})$ to an $h$-space contradicts evasivity.
\end{proof}

\section{Auxiliary fields and support inequalities}\label{sec:auxiliary}

We use the conventions of Section~\ref{sec:preliminaries}.
Choose a positive integer $s$ relatively prime to $n$, and regard
$\mathbb F_{q^n}$ and $\mathbb F_{q^s}$ as subfields of
$\mathbb F_{q^{ns}}$.
The multiplication map identifies $\mathbb F_{q^n}\otimes_{\mathbb F_q}\mathbb F_{q^s}$ with $\mathbb F_{q^{ns}}$.
In particular, $\mathbb F_{q^n}\cap \mathbb F_{q^s}=\mathbb F_q$, and a basis of $\mathbb F_{q^n}/\mathbb F_q$ is
also a basis of $\mathbb F_{q^{ns}}/\mathbb F_{q^s}$. If $B\le_{\mathbb F_q}\mathbb F_{q^s}^T$, this identification
gives an injection
\[
 \mathbb F_{q^n}\otimes_{\mathbb F_q}B\longrightarrow \mathbb F_{q^{ns}}^T,
 \qquad e\otimes b\longmapsto eb.
\]
Consequently $\dim_{\mathbb F_{q^n}}\langle B\rangle_{\mathbb F_{q^n}}=\dim_{\mathbb F_q} B$. We use this observation
both for relation spaces and for their supports.

The scalar-orbit mechanism in the case $h=k\le n-2$ is developed in
Giannoni and Neri~\cite{GN}. Their finite \'etale algebra argument also
removes the coprimality restriction on the auxiliary degree in that
case. Here we keep $\gcd(n,s)=1$ and control relation spaces of higher
dimension. The following intersection lemma supplies the required
reduction to relations with small support.

\begin{lemma}\label{lem:intersection}
Let $V=\mathbb F_{q^{ns}}^y$ and let $E\le_{\mathbb F_{q^n}} V$, with $\dim_{\mathbb F_{q^n}} E=x$. There exist
$\mathbb F_{q^{ns}}$-subspaces $P_0,\ldots,P_{s-1}\le_{\mathbb F_{q^{ns}}} V$, depending only on $E$, such
that
\[
 \sum_{v=0}^{s-1}\dim_{\mathbb F_{q^{ns}}} P_v=x,
 \qquad
 \dim_{\mathbb F_q}(E\cap\Gamma)
 \le\sum_{v=0}^{s-1}\dim_{\mathbb F_{q^s}}(\Gamma\cap P_v)
\]
for every $\mathbb F_{q^s}$-subspace $\Gamma\le_{\mathbb F_{q^s}} V$.
\end{lemma}

\begin{proof}
Give $V\otimes_{\mathbb F_q}\mathbb F_{q^s}$ its $\mathbb F_{q^{ns}}$-vector space structure by letting
$\mathbb F_{q^n}$ act on the first factor and $\mathbb F_{q^s}$ on the second. Let
$\phi_0,\ldots,\phi_{s-1}$ be the elements of $\operatorname{Gal}(\mathbb F_{q^{ns}}/\mathbb F_{q^n})$.
Their restrictions to $\mathbb F_{q^s}$ are the elements of
$\operatorname{Gal}(\mathbb F_{q^s}/\mathbb F_q)$. The separable splitting of
$\mathbb F_{q^s}\otimes_{\mathbb F_q}\mathbb F_{q^s}$ yields the $\mathbb F_{q^{ns}}$-linear isomorphism
\[
 \Psi:V\otimes_{\mathbb F_q}\mathbb F_{q^s}\longrightarrow\bigoplus_{v=0}^{s-1}V,
 \qquad
 z\otimes\lambda\longmapsto
 \bigl(\phi_v^{-1}(z)\lambda\bigr)_{v=0}^{s-1},
\]
where the automorphisms act coordinatewise. Composing its $v$th
coordinate with $\phi_v$ gives $z\otimes\lambda\mapsto
(z\phi_v(\lambda))_v$, the standard Galois splitting, so $\Psi$ is
bijective. Since every $\phi_v$ fixes $\mathbb F_{q^n}$, the two
specified scalar actions also give its asserted
$\mathbb F_{q^{ns}}$-linearity.

Put
\[
 B=\Psi(E\otimes_{\mathbb F_q}\mathbb F_{q^s}),\qquad
 \Gamma_v=\phi_v^{-1}(\Gamma).
\]
Then $B$ is an $\mathbb F_{q^{ns}}$-subspace of dimension $x$. The same splitting, applied
to an $\mathbb F_{q^s}$-basis of $\Gamma$, gives the equality of $\mathbb F_{q^s}$-subspaces
\[
 \Psi(\Gamma\otimes_{\mathbb F_q}\mathbb F_{q^s})=\bigoplus_{v=0}^{s-1}\Gamma_v.
\]
Extension of scalars preserves the intersection of two vector
subspaces. Hence
\[
 I:=\Psi((E\cap\Gamma)\otimes_{\mathbb F_q}\mathbb F_{q^s})
   =B\cap\bigoplus_{v=0}^{s-1}\Gamma_v,
 \qquad \dim_{\mathbb F_{q^s}} I=\dim_{\mathbb F_q}(E\cap\Gamma).
\]

For $0\le v\le s$, define
\[
 B_v=\{b\in B:b_0=\cdots=b_{v-1}=0\}.
\]
Thus $B_0=B$ and $B_s=0$. If $\pi_v$ is the $v$th coordinate projection,
set $Q_v=\pi_v(B_v)$ for $v<s$. Since the kernel of
$\pi_v|_{B_v}$ is $B_{v+1}$,
\[
 \sum_{v=0}^{s-1}\dim_{\mathbb F_{q^{ns}}} Q_v=\dim_{\mathbb F_{q^{ns}}} B=x.
\]
With $I_v=I\cap B_v$, the kernel of $\pi_v|_{I_v}$ is $I_{v+1}$ and
$\pi_v(I_v)\subseteq Q_v\cap\Gamma_v$. Therefore
\[
 \dim_{\mathbb F_q}(E\cap\Gamma)
 =\sum_{v=0}^{s-1}\dim_{\mathbb F_{q^s}}\pi_v(I_v)
 \le\sum_{v=0}^{s-1}\dim_{\mathbb F_{q^s}}(Q_v\cap\Gamma_v).
\]
Take $P_v=\phi_v(Q_v)$. Semilinear application of $\phi_v$ preserves
the relevant dimensions and changes $\Gamma_v$ into $\Gamma$.
The construction of $B_v$, $Q_v$ and $P_v$ uses only $E$.
\end{proof}

\begin{remark}\label{rem:scalar-intersection}
The case $y=1$ of Lemma~\ref{lem:intersection} already gives a product
bound. In that case every $P_v$ is either $0$ or $\mathbb F_{q^{ns}}$,
and exactly $\dim_{\mathbb F_{q^n}} E$ of them are nonzero. Hence, for
$E\le_{\mathbb F_{q^n}} \mathbb F_{q^{ns}}$ and
$\Gamma\le_{\mathbb F_{q^s}} \mathbb F_{q^{ns}}$,
\[
 \dim_{\mathbb F_q}(E\cap\Gamma)\le (\dim_{\mathbb F_{q^n}} E)(\dim_{\mathbb F_{q^s}}\Gamma).
\]
\end{remark}

Fix a basis $e_1,\ldots,e_n$ of $\mathbb F_{q^n}/\mathbb F_q$. Every vector $w\in \mathbb F_{q^{ns}}^T$
has a unique expansion $w=\sum_{b=1}^n e_b w_b$ with $w_b\in \mathbb F_{q^s}^T$.
For an $\mathbb F_{q^n}$-subspace $W\le_{\mathbb F_{q^n}} \mathbb F_{q^{ns}}^T$, define
\[
 \operatorname{Supp}_q(W)
 =\langle w_b:w\in W,\ 1\le b\le n\rangle_{\mathbb F_q},
 \qquad
 \sigma_q(W)=\dim_{\mathbb F_q}\operatorname{Supp}_q(W).
\]
For an $\mathbb F_{q^{ns}}$-subspace $Z\le_{\mathbb F_{q^{ns}}} \mathbb F_{q^{ns}}^T$, define instead
\[
 \operatorname{Supp}_{q^s}(Z)
 =\langle z_b:z\in Z,\ 1\le b\le n\rangle_{\mathbb F_{q^s}},
 \qquad
 \widetilde\sigma(Z)=\dim_{\mathbb F_{q^s}}\operatorname{Supp}_{q^s}(Z).
\]
A change of the basis $(e_b)$ multiplies the list of components by an
invertible matrix over $\mathbb F_q$, so both definitions are independent
of that choice. Moreover,
\[
 W\subseteq\langle\operatorname{Supp}_q(W)\rangle_{\mathbb F_{q^n}},
 \qquad
 Z\subseteq\langle\operatorname{Supp}_{q^s}(Z)\rangle_{\mathbb F_{q^{ns}}}.
\]
In particular, $\widetilde\sigma(Z)\ge\dim_{\mathbb F_{q^{ns}}} Z$.

\begin{lemma}\label{lem:density}
Let $W\le_{\mathbb F_{q^n}} \mathbb F_{q^{ns}}^T$, write $x=\dim_{\mathbb F_{q^n}} W$, and put $Y=\langle W\rangle_{\mathbb F_{q^{ns}}}$.
There are $\mathbb F_{q^{ns}}$-subspaces $Z_0,\ldots,Z_{s-1}\le_{\mathbb F_{q^{ns}}} Y$ such that
\[
 \sum_{v=0}^{s-1}\dim_{\mathbb F_{q^{ns}}} Z_v=x,
 \qquad
 \sigma_q(W)\ge\sum_{v=0}^{s-1}\widetilde\sigma(Z_v).
\]
If $W\ne0$, it follows that some nonzero $\mathbb F_{q^{ns}}$-subspace $Z\le_{\mathbb F_{q^{ns}}} Y$ satisfies
\[
 \frac{\widetilde\sigma(Z)}{\dim_{\mathbb F_{q^{ns}}} Z}
 \le\frac{\sigma_q(W)}{\dim_{\mathbb F_{q^n}} W}.
\]
\end{lemma}

\begin{proof}
The assertion for $W=0$ follows by taking every $Z_v=0$. Suppose
$W\ne0$, set $y=\dim_{\mathbb F_{q^{ns}}} Y$, and let the rows of $M\in \mathbb F_{q^{ns}}^{y\times T}$ be
an $\mathbb F_{q^{ns}}$-basis of $Y$. The map $\xi\mapsto\xi M$ identifies an $\mathbb F_{q^n}$-subspace
$E\le_{\mathbb F_{q^n}} \mathbb F_{q^{ns}}^y$ of dimension $x$ with $W$. Put
\[
 \Gamma=M\mathbb F_{q^s}^T\le_{\mathbb F_{q^s}} \mathbb F_{q^{ns}}^y,\qquad c=\dim_{\mathbb F_{q^s}}\Gamma.
\]
On $\mathbb F_{q^{ns}}^y$, use the nondegenerate $\mathbb F_{q^n}$-bilinear pairing
\[
 (\xi,v)\longmapsto\operatorname{Tr}_{\mathbb F_{q^{ns}}/\mathbb F_{q^n}}(\xi\cdot v),
\]
and let $E^\perp$ be the annihilator of $E$. Thus
$\dim_{\mathbb F_{q^n}} E^\perp=sy-x$.

Set $S=\operatorname{Supp}_q(W)$, and take its annihilator in $\mathbb F_{q^s}^T$
with respect to the nondegenerate $\mathbb F_q$-bilinear trace pairing.
For $e\in \mathbb F_{q^n}$ and $\lambda\in \mathbb F_{q^s}$ one has
\[
 \operatorname{Tr}_{\mathbb F_{q^{ns}}/\mathbb F_{q^n}}(e\lambda)
   =e\operatorname{Tr}_{\mathbb F_{q^s}/\mathbb F_q}(\lambda).
\]
Indeed, the restrictions of $\operatorname{Gal}(\mathbb F_{q^{ns}}/\mathbb F_{q^n})$ to $\mathbb F_{q^s}$ exhaust
$\operatorname{Gal}(\mathbb F_{q^s}/\mathbb F_q)$. Expanding a vector $w$ in the basis
$(e_b)$ now gives
\[
 S^\perp
 =\{u\in \mathbb F_{q^s}^T:\operatorname{Tr}_{\mathbb F_{q^{ns}}/\mathbb F_{q^n}}(w\cdot u)=0
                    \text{ for every }w\in W\}.
\]
This is the inverse image of $\Gamma\cap E^\perp$ under $u\mapsto Mu$.
That map has an $\mathbb F_{q^s}$-kernel of dimension $T-c$, whence
\[
 \sigma_q(W)=sc-\dim_{\mathbb F_q}(\Gamma\cap E^\perp).
\]
Apply Lemma~\ref{lem:intersection} to $E^\perp$. Its subspaces $P_v$
satisfy $\sum_v\dim_{\mathbb F_{q^{ns}}} P_v=sy-x$, and therefore
\[
 \sigma_q(W)\ge
 \sum_{v=0}^{s-1}\bigl(c-\dim_{\mathbb F_{q^s}}(\Gamma\cap P_v)\bigr).
\]
Let $R_v$ be the annihilator of $P_v$ for the ordinary $\mathbb F_{q^{ns}}$-bilinear dot
product, and set $Z_v=R_vM$. Since multiplication by $M$ is injective,
\[
 \dim_{\mathbb F_{q^{ns}}} Z_v=y-\dim_{\mathbb F_{q^{ns}}} P_v,
 \qquad\sum_v\dim_{\mathbb F_{q^{ns}}} Z_v=x.
\]

To identify the remaining summands, let $\mathbf R_v$ have an $\mathbb F_{q^{ns}}$-basis of
$R_v$ as its rows. Expand $\mathbf R_vM$ over the $\mathbb F_{q^s}$-basis $(e_b)$ of $\mathbb F_{q^{ns}}$.
The $\mathbb F_{q^s}$-row span of the expanded matrix is
$\operatorname{Supp}_{q^s}(Z_v)$: multiplication of a row by a scalar in
$\mathbb F_{q^{ns}}$ only takes $\mathbb F_{q^s}$-linear combinations of its expanded rows. The same
expanded matrix represents the $\mathbb F_{q^s}$-linear map
$u\mapsto \mathbf R_vMu$ from $\mathbb F_{q^s}^T$. Row rank equals column rank, so
\[
 \widetilde\sigma(Z_v)
 =\dim_{\mathbb F_{q^s}} \mathbf R_v(\Gamma)
 =c-\dim_{\mathbb F_{q^s}}(\Gamma\cap P_v).
\]
This proves the first assertion. Since $x>0$, at least one $Z_v$ is
nonzero. If all the nonzero $Z_v$ had density strictly greater than
$\sigma_q(W)/x$, summing their support dimensions would contradict
the first assertion. This proves the last claim.
\end{proof}

\section{A finite criterion for evasive systems}\label{sec:finite-evasive}

Fix integers $1\le h\le k$, and write
\[
 a=k+1,\qquad j=a-h,\qquad \alpha=\frac ja.
\]
The coprimality hypothesis $\gcd(h,k+1)=1$ is equivalent to
$\gcd(a,j)=1$. The construction below requires $nj\ge a$; it does
not require $k<n$. Under coprimality, $nj=a$ holds if and only if
$j=1$ and $n=a$. Choose $\ell\ge1$ with
\[
 T=n\ell\alpha\in\mathbb Z,\qquad r=s\ell,\qquad N=sT.
\]
For a matrix $A\in \mathbb F_{q^{ns}}^{\ell\times T}$, considered as a map on column
vectors, set
\[
 \mathcal K=\ker_{\mathbb F_{q^{ns}}} A,\qquad U=A\mathbb F_{q^s}^T\le_{\mathbb F_{q^s}} \mathbb F_{q^{ns}}^\ell.
\]
The ambient space $\mathbb F_{q^{ns}}^\ell$ has dimension $r$ over $\mathbb F_{q^n}$.

\begin{lemma}\label{lem:kernel-evasivity}
Suppose $sT\ge a$, $s\ell\ge h+1$, and $A|_{\mathbb F_{q^s}^T}$ is injective. Then
$U=A\mathbb F_{q^s}^T$ is $(h,k)$-evasive if and only if there is no $\mathbb F_{q^n}$-subspace
$W\le_{\mathbb F_{q^n}} \mathcal K$ with
\[
 \dim_{\mathbb F_{q^n}} W=j,\qquad\sigma_q(W)\le a.
\]
\end{lemma}

\begin{proof}
Suppose an $\mathbb F_{q^n}$-subspace of dimension $h$ contains $a$ vectors of $U$
independent over $\mathbb F_q$. Their preimages span a subspace
$B\le_{\mathbb F_q}\mathbb F_{q^s}^T$ of dimension $a$. Linear disjointness gives
$\dim_{\mathbb F_{q^n}}\langle B\rangle_{\mathbb F_{q^n}}=a$, whereas the image of this $\mathbb F_{q^n}$-space under
$A$ has dimension at most $h$. Its kernel therefore contains a
$j$-dimensional $\mathbb F_{q^n}$-subspace $W$. Expanding elements of
$\langle B\rangle_{\mathbb F_{q^n}}$ over $(e_b)$ shows that
$\operatorname{Supp}_q(W)\subseteq B$, so $\sigma_q(W)\le a$.

Conversely, let $W$ satisfy the displayed conditions. Extend
$\operatorname{Supp}_q(W)$ to an $a$-dimensional
$\mathbb F_q$-subspace $B\le \mathbb F_{q^s}^T$, which is possible because $sT\ge a$.
Then $W\subseteq\langle B\rangle_{\mathbb F_{q^n}}$, and the dimension formula gives
$\dim_{\mathbb F_{q^n}} A\langle B\rangle_{\mathbb F_{q^n}}\le a-j=h$. The subspace $A(B)$ has
$\mathbb F_q$-dimension $a$, by injectivity. Extending its $\mathbb F_{q^n}$-span to
an $h$-dimensional ambient subspace contradicts $(h,k)$-evasivity.
\end{proof}

For later use define the finite index set
\[
 \mathcal I_0=
 \{(d,\lfloor ad/j\rfloor):1\le d<j\}.
\]
It is empty if $j=1$. Only the largest permissible support dimension
for each $d$ is needed, because the count below covers all smaller
supports. Also put $\Delta=nhj-a^2$, and define
\[
 \begin{aligned}
 C_b(Q)&=\prod_{u=1}^b(1-Q^{-u})^{-1}\quad(b\ge0,\ Q>1),\\
 E_\ell(d,c)&=n\ell(\alpha c-d)+nd(c-d)-c^2,
 \end{aligned}
\]
with $C_0(Q)=1$.

\begin{theorem}\label{thm:finite-evasive}
Assume $1\le h\le k$, $\gcd(a,j)=1$, $nj\ge a$, and
$\gcd(n,s)=1$. Suppose $s\ell\ge h+1$ and
$T=n\ell j/a\in\mathbb Z$ with $sT\ge a$ and $T\ge c$ for every
$(d,c)\in\mathcal I_0$.
Define
\[
 \begin{aligned}
 G_{q,s,\ell}
   &=\gauss{sT}{a}{q}\gauss{a}{j}{q^n}
       \frac{q-1}{q^s-1}\,q^{-ns\ell j},\\
 D_{q,s,\ell}(d,c)
   &=\gauss{T}{c}{q^s}\gauss{c}{d}{q^{ns}}
       q^{-ns\ell d},\\
 R_{q,s,\ell}
   &=\begin{cases}
       0,&n=a\text{ and }j=1,\\
       \displaystyle\frac{q^{ns\ell}-1}{q^{ns}-1}\,q^{-nsT},
          &\text{otherwise}.
      \end{cases}
 \end{aligned}
\]
If
\[
 \mathcal B_{q,s,\ell}:=
 G_{q,s,\ell}
 +\sum_{(d,c)\in\mathcal I_0}D_{q,s,\ell}(d,c)
 +R_{q,s,\ell}<1,
\]
there is an $\mathbb F_{q^s}$-linear subspace $U\le \mathbb F_{q^{ns}}^\ell$ which spans $\mathbb F_{q^{ns}}^\ell$
over $\mathbb F_{q^n}$, is $(h,k)$-evasive, and has
\[
 \dim_{\mathbb F_q} U=N=\frac{rnj}{a},\qquad r=s\ell.
\]
If $k<n$, this system is maximum. The following finite-product
bounds will be useful:
\begin{equation}\label{eq:refined-counts}
 \begin{aligned}
 G_{q,s,\ell}
 &\le \frac{q^{\Delta+1-s}}{1-q^{-s}}
       \prod_{u=2}^a(1-q^{-u})^{-1}C_j(q^n),\\
 D_{q,s,\ell}(d,c)
 &\le q^{sE_\ell(d,c)}C_c(q^s)C_d(q^{ns}).
 \end{aligned}
\end{equation}
Each right-hand side is nonincreasing in $q\ge2$ whenever its
exponent of $q$ is nonpositive. A simpler sufficient condition is
\[
 16q^{\Delta-s+1}
 +16\sum_{(d,c)\in\mathcal I_0}
    q^{s[n\ell(\alpha c-d)+nd(c-d)-c^2]}
 +R_{q,s,\ell}<1.
\]
For fixed $q,n,h,k,s$, the condition $16q^{\Delta-s+1}<1$ consequently
gives such a system for every sufficiently large $\ell$ with
$a\mid n\ell j$.
\end{theorem}

\begin{proof}
Choose $A$ uniformly from $\mathbb F_{q^{ns}}^{\ell\times T}$. For a fixed $\mathbb F_{q^{ns}}$-subspace
$X\le_{\mathbb F_{q^{ns}}} \mathbb F_{q^{ns}}^T$ of dimension $d$, independence of the rows of $A$ gives
\[
 \Pr(X\subseteq \mathcal K)=q^{-ns\ell d}.
\]
Call an $\mathbb F_{q^n}$-subspace $W$ a bad relation space if
$\dim_{\mathbb F_{q^n}} W=j$ and $\sigma_q(W)\le a$. Separate these spaces according
to whether $\dim_{\mathbb F_{q^{ns}}}\langle W\rangle_{\mathbb F_{q^{ns}}}=j$.

First suppose $W\le \mathcal K$ is bad and
$y=\dim_{\mathbb F_{q^{ns}}}\langle W\rangle_{\mathbb F_{q^{ns}}}<j$. Lemma~\ref{lem:density} supplies a
nonzero $\mathbb F_{q^{ns}}$-space $X\le\langle W\rangle_{\mathbb F_{q^{ns}}}\le \mathcal K$ for which, on putting
$d=\dim_{\mathbb F_{q^{ns}}}X$ and $c=\lfloor ad/j\rfloor$,
\[
 1\le d<j,\qquad d\le\widetilde\sigma(X)\le c,
 \qquad(d,c)\in\mathcal I_0.
\]
The number of $d$-dimensional $\mathbb F_{q^{ns}}$-spaces with support dimension at
most $c$ is at most
\[
 \gauss{T}{c}{q^s}\gauss{c}{d}{q^{ns}}.
\]
Indeed, choose a $c$-dimensional $\mathbb F_{q^s}$-space $S\le \mathbb F_{q^s}^T$ containing the
support, and then choose $X\le_{\mathbb F_{q^{ns}}}\langle S\rangle_{\mathbb F_{q^{ns}}}$, whose ambient
$\mathbb F_{q^{ns}}$-dimension is $c$. This may count a space several times, which is
harmless for an upper bound. Since $T\ge c$, every required
extension of a support exists. A union bound proves that the
probability of such an $X\le \mathcal K$ is at most $D_{q,s,\ell}(d,c)$.

Now consider the bad spaces $W$ with
$\dim_{\mathbb F_{q^{ns}}}\langle W\rangle_{\mathbb F_{q^{ns}}}=j$. Their total number is at most
\[
 \gauss{sT}{a}{q}\gauss{a}{j}{q^n}:
\]
choose an $a$-dimensional $\mathbb F_q$-space $B\le \mathbb F_{q^s}^T$ containing
the support, which is possible because $sT\ge a$, then a $j$-dimensional $\mathbb F_{q^n}$-space in
$\langle B\rangle_{\mathbb F_{q^n}}$, which has $\mathbb F_{q^n}$-dimension $a$.

Multiplication by $\mathbb F_{q^s}^*$ preserves this class of spaces and their
support dimensions. Its stabilizer at every such $W$ is exactly
$\mathbb F_q^*$. To see this, let $w_1,\ldots,w_j$ be an $\mathbb F_{q^n}$-basis of
$W$. It is also $\mathbb F_{q^{ns}}$-linearly independent. If $\mu\in \mathbb F_{q^s}^*$ and
$\mu W=W$, write $\mu w_1=\sum_i b_iw_i$ with $b_i\in \mathbb F_{q^n}$.
Independence over $\mathbb F_{q^{ns}}$ forces $\mu=b_1\in \mathbb F_{q^n}$, so
$\mu\in \mathbb F_{q^n}\cap \mathbb F_{q^s}=\mathbb F_q$; the converse is immediate.
Thus every orbit has size $(q^s-1)/(q-1)$. This scalar-orbit
count is the mechanism used for scattered subspaces by Giannoni
and Neri~\cite{GN}. The events
$W\subseteq \mathcal K$ are constant on these orbits because $\mathcal K$ is $\mathbb F_{q^{ns}}$-linear,
and each has probability $q^{-ns\ell j}$. Counting one representative
per orbit gives the upper bound $G_{q,s,\ell}$ for this class of
failures. The orbit division has only been used for spaces whose
$\mathbb F_{q^n}$-dimension equals their $\mathbb F_{q^{ns}}$-span dimension.

Avoiding these relation events also makes $A|_{\mathbb F_{q^s}^T}$ injective.
Indeed, a nonzero $v\in \mathcal K\cap \mathbb F_{q^s}^T$ gives an $\mathbb F_{q^{ns}}$-line $\mathbb F_{q^{ns}}v$ of support
dimension one. If $j>1$, it is included in the event indexed by
$(1,\lfloor a/j\rfloor)\in\mathcal I_0$, because $\lfloor a/j\rfloor\ge1$.
If $j=1$, the $\mathbb F_{q^n}$-line $\mathbb F_{q^n}v$ is a bad relation space in the second
class. Lemma~\ref{lem:kernel-evasivity} then gives evasivity, and
injectivity gives $\dim_{\mathbb F_q} U=sT$.

It remains to ensure full span. Since $nj\ge a$, we have
$T\ge\ell$. Equality means $nj=a$, which forces $j=1$ and $n=a$
because $\gcd(a,j)=1$. If
$T>\ell$, a failure of full row rank gives a line in the left kernel
of $A$. There are $(q^{ns\ell}-1)/(q^{ns}-1)$ possible lines, and a
fixed one annihilates all $T$ columns with probability $q^{-nsT}$.
This is the stated bound $R_{q,s,\ell}$. In the endpoint $T=\ell$,
any nonzero $w\in \mathcal K$ has at most $n=a$ components in its expansion
over $(e_b)$. The support of $\mathbb F_{q^n}w$ is contained in their
$\mathbb F_q$-span, because multiplication by an element of $\mathbb F_{q^n}$
acts on the components through a matrix over $\mathbb F_q$.
Thus $\mathbb F_{q^n}w$ is already a bad relation space in the second class.
Avoiding that class gives $\mathcal K=0$ and full row rank without a separate
failure event, explaining the value $R_{q,s,\ell}=0$ at the endpoint.

Full row rank now implies full $\mathbb F_{q^n}$-span for the particular subspace
$U=A\mathbb F_{q^s}^T$. Explicitly, $\langle \mathbb F_{q^s}\rangle_{\mathbb F_{q^n}}=\mathbb F_{q^{ns}}$, so
$\langle \mathbb F_{q^s}^T\rangle_{\mathbb F_{q^n}}=\mathbb F_{q^{ns}}^T$, and $A$ is $\mathbb F_{q^n}$-linear. Consequently
\[
 \langle U\rangle_{\mathbb F_{q^n}}
 =A\langle \mathbb F_{q^s}^T\rangle_{\mathbb F_{q^n}}=A(\mathbb F_{q^{ns}}^T)=\mathbb F_{q^{ns}}^\ell.
\]
This proves the assertion whenever $\mathcal B_{q,s,\ell}<1$.
When $k<n$, the dimension equals the upper bound in
Theorem~\ref{thm:evasive-bound}.

It remains to prove \eqref{eq:refined-counts} and the remaining
assertions of the theorem. For every prime power $Q$ and
$0\le v\le u$,
\[
 \gauss{u}{v}{Q}
 \le Q^{v(u-v)}C_v(Q)
 \le Q^{v(u-v)}\prod_{i\ge1}(1-Q^{-i})^{-1}
 <4Q^{v(u-v)}.
\]
The last estimate follows, for example, from
\[
 \prod_{i\ge1}(1-2^{-i})
 \ge(1-2^{-1})(1-2^{-2})
       \left(1-\sum_{i\ge3}2^{-i}\right)
 =\frac9{32}>\frac14.
\]
The product estimate follows first for finite products from
$\prod_i(1-x_i)\ge1-\sum_i x_i$, and then by passage to the limit.
In the bound for $G$, the factor $1-q^{-1}$ from
$(q-1)/(q^s-1)=q^{1-s}(1-q^{-1})/(1-q^{-s})$ cancels the
first factor of $C_a(q)$. This proves \eqref{eq:refined-counts},
using $aT=n\ell j$; the assertion about monotonicity follows
factor by factor. Using also $(q-1)/(q^s-1)\le q^{1-s}$ gives
\[
 G_{q,s,\ell}<16q^{nhj-a^2-s+1},
 \qquad
 D_{q,s,\ell}(d,c)
 <16q^{s[n\ell(\alpha c-d)+nd(c-d)-c^2]}.
\]
For $(d,c)\in\mathcal I_0$, coprimality gives $jc<ad$: equality
would imply $j\mid d$, contrary to $0<d<j$. Every exponent in the
second bound therefore tends to $-\infty$ with $\ell$. When
$T>\ell$, the remaining term satisfies
\[
 R_{q,s,\ell}
 <\frac{q^{-ns(T-\ell)}}{q^{ns}-1}\longrightarrow0.
\]
There are finitely many pairs $(d,c)$, so a fixed margin below one
in the generic bound proves the final assertion.
\end{proof}

\section{Simultaneous constraints and a uniform threshold}\label{sec:simultaneous}

The density reduction can impose additional constraints whose ratios
are strictly smaller than the one determining the dimension. It does
not require coprimality for those additional constraints.

\begin{theorem}\label{thm:simultaneous}
Let $(h_i,k_i)$, $0\le i\le \kappa$, be a finite list with
$1\le h_i\le k_i$. Set
\[
 a_i=k_i+1,\qquad j_i=a_i-h_i,\qquad \rho_i=\frac{a_i}{j_i},
 \qquad \Delta_i=nh_ij_i-a_i^2.
\]
Assume
\[
 \gcd(a_0,j_0)=1,\qquad nj_0\ge a_0,\qquad
 \rho_i<\rho_0\quad(1\le i\le \kappa).
\]
Write $a=a_0$, $j=j_0$, $\alpha=j/a$ and $T=n\ell\alpha$.
Let $s$ be coprime to $n$, assume
$s\ell\ge1+\max_i h_i$, and suppose $T$ is an integer with
$sT\ge\max_i a_i$. For $1\le i\le \kappa$, put
\[
 \mathcal I_i=
 \{(d,\lfloor a_i d/j_i\rfloor):1\le d< j_i\},
 \qquad
 G^{(i)}_{q,s,\ell}
   =\gauss{sT}{a_i}{q}\gauss{a_i}{j_i}{q^n}
       \frac{q-1}{q^s-1}\,q^{-ns\ell j_i},
\]
and let $F^{(i)}_{q,s,\ell}=G^{(i)}_{q,s,\ell}$ when $T<a_i$, and
$F^{(i)}_{q,s,\ell}=\min\{G^{(i)}_{q,s,\ell},D_{q,s,\ell}(j_i,a_i)\}$
when $T\ge a_i$. Assume that $T\ge c$ for every pair $(d,c)$ in
$\mathcal I_0\cup\cdots\cup\mathcal I_\kappa$.
Use $\mathcal B_{q,s,\ell}$ and $D_{q,s,\ell}(d,c)$ from
Theorem~\ref{thm:finite-evasive}, with the parameters of index $0$.
If
\[
 \mathcal B_{q,s,\ell}
 +\sum_{i=1}^{\kappa}\Bigl(F^{(i)}_{q,s,\ell}
   +\sum_{(d,c)\in\mathcal I_i}D_{q,s,\ell}(d,c)\Bigr)<1,
\]
there is an $\mathbb F_{q^s}$-linear subspace $U\le \mathbb F_{q^{ns}}^\ell$ spanning $\mathbb F_{q^{ns}}^\ell$
over $\mathbb F_{q^n}$, of dimension
\[
 \dim_{\mathbb F_q} U=\frac{sn\ell j}{a},
\]
which is $(h_i,k_i)$-evasive for every $i$. Moreover,
\begin{equation}\label{eq:secondary-generic}
 G^{(i)}_{q,s,\ell}
 \le\frac{q^{\Delta_i+1-s-ns\ell(j_i-\alpha a_i)}}{1-q^{-s}}
       \prod_{u=2}^{a_i}(1-q^{-u})^{-1}C_{j_i}(q^n),
\end{equation}
where $j_i-\alpha a_i=j_i(1-\rho_i/\rho_0)>0$. In particular, if
$16q^{nh_0j-a^2-s+1}<1$, the conclusion holds for every sufficiently
large $\ell$ for which $a\mid n\ell j$.
\end{theorem}

\begin{proof}
Use the same random matrix $A$ as in
Theorem~\ref{thm:finite-evasive}; the hypotheses $sT\ge a$ and $T\ge c$
for $(d,c)\in\mathcal I_0$ are those used there. Its primary relation
events ensure injectivity on $\mathbb F_{q^s}^T$ and full span. By
Lemma~\ref{lem:kernel-evasivity}, which only uses $sT\ge a_i$ and
$s\ell\ge h_i+1$, failure of the constraint of index
$i\ge1$ produces an $\mathbb F_{q^n}$-subspace $W\le \mathcal K$ with
$\dim_{\mathbb F_{q^n}} W=j_i$ and $\sigma_q(W)\le a_i$. As for the
primary constraint, we separate these spaces according to the
dimension of their span over $\mathbb F_{q^{ns}}$.

Suppose first that $\dim_{\mathbb F_{q^{ns}}}\langle W\rangle_{\mathbb F_{q^{ns}}}=j_i$.
The count in the proof of Theorem~\ref{thm:finite-evasive} applies with
$(a,j)$ replaced by $(a_i,j_i)$. Since $sT\ge a_i$, there are at most
$\gauss{sT}{a_i}{q}\gauss{a_i}{j_i}{q^n}$ such spaces. Multiplication
by $\mathbb F_{q^s}^*$ preserves this class, and its stabilizer at every
member is exactly $\mathbb F_q^*$: the argument given there uses only the
equality between the $\mathbb F_{q^n}$-dimension of $W$ and the dimension
of its $\mathbb F_{q^{ns}}$-span. Each event $W\subseteq\mathcal K$ has
probability $q^{-ns\ell j_i}$ and is constant on orbits, which gives the
bound $G^{(i)}_{q,s,\ell}$. When $T\ge a_i$, the space $Z=\langle W\rangle_{\mathbb F_{q^{ns}}}$
satisfies $\operatorname{Supp}_{q^s}(Z)
=\langle\operatorname{Supp}_q(W)\rangle_{\mathbb F_{q^s}}$, so $Z$ is a
$j_i$-dimensional subspace of $\mathcal K$ with
$\widetilde\sigma(Z)\le a_i$. This gives the alternative bound
$D_{q,s,\ell}(j_i,a_i)$ for the same class.

Suppose next that $\dim_{\mathbb F_{q^{ns}}}\langle W\rangle_{\mathbb F_{q^{ns}}}<j_i$.
Lemma~\ref{lem:density} gives an $\mathbb F_{q^{ns}}$-subspace $X\le \mathcal K$ with
\[
 1\le d:=\dim_{\mathbb F_{q^{ns}}}X< j_i,\qquad
 d\le\widetilde\sigma(X)\le c:=\lfloor a_i d/j_i\rfloor.
\]
Since $T\ge c$, the count of supports of dimension at most $c$ bounds
its occurrence by $D_{q,s,\ell}(d,c)$. The crucial coefficient is now
\[
 \alpha c-d\le
 d\left(\frac{\rho_i}{\rho_0}-1\right)<0.
\]

The estimate \eqref{eq:secondary-generic} follows from the same
finite-product bounds as \eqref{eq:refined-counts}, together with
$a_i(sT-a_i)+nj_ih_i-ns\ell j_i=\Delta_i-ns\ell(j_i-\alpha a_i)$, which
holds because $aT=n\ell j$. Thus every additional term tends to zero
with $\ell$. Adding the union bounds proves both claims.
\end{proof}

The following explicit threshold will also make the later passage to
all admissible ranks effective. It gives a computable bound independent of $q$.

\begin{proposition}\label{prop:uniform-threshold}
Retain the finite list and assumptions of
Theorem~\ref{thm:simultaneous}, allowing $\kappa=0$. Put
\[
 \lambda=\frac{nj}{a},\qquad
 A_* =\max_i a_i,\qquad H_* =\max_i h_i,\qquad
 \Delta=nh_0j-a^2.
\]
Choose $s$ with $\gcd(n,s)=1$ and $s\ge\Delta+6$. Let
$\mathcal J$ be the disjoint union of the sets $\mathcal I_i$ for
$0\le i\le \kappa$; repetitions of a pair in different sets are retained.
Thus $|\mathcal J|=\sum_{i=0}^{\kappa}(j_i-1)$.
For integers $d,c$ set
\[
 \varepsilon_{d,c}=d-\frac ja c,
 \qquad B_{d,c}=nd(c-d)-c^2,
 \qquad \mathcal M=|\mathcal J|+\kappa+\boldsymbol{1}_{\lambda>1}.
\]
Start with
\[
 L_s=\max\left\{1,
          \left\lceil\frac{A_*}{\lambda}\right\rceil,
          \left\lceil\frac{H_*+1}{s}\right\rceil\right\}.
\]
If $\mathcal M>0$, enlarge this maximum by including, for every pair in
$\mathcal J$, the integer
\[
 \left\lceil
 \frac{B_{d,c}+s^{-1}\log_2(64\mathcal M)}{n\varepsilon_{d,c}}
 \right\rceil,
\]
for every $1\le i\le\kappa$ the integer
\[
 \left\lceil
 \frac{\Delta_i+1-s+\log_2(64\mathcal M)}{ns\,\varepsilon_{j_i,a_i}}
 \right\rceil,
\]
and, if $\lambda>1$, also the integer
\[
 \left\lceil\frac{\log_2(4\mathcal M)}{ns(\lambda-1)}\right\rceil.
\]
Then for every prime power $q$ and every integer $\ell\ge L_s$ with
$a\mid n\ell j$, the conclusion of
Theorem~\ref{thm:simultaneous} holds. Thus $L_s$ depends only on the
listed integer parameters, not on $q$.
\end{proposition}

\begin{proof}
The initial entries in the maximum give $T\ge A_*$, hence $sT\ge A_*$
and $T\ge c$ for every pair in $\mathcal J$, together with
$s\ell\ge H_*+1$. The positivity of $\varepsilon_{d,c}$ for
$(d,c)\in\mathcal J$ was proved in Theorem~\ref{thm:finite-evasive} for
$i=0$ and in Theorem~\ref{thm:simultaneous} for $i\ge1$, which also gives
$\varepsilon_{j_i,a_i}=j_i-\alpha a_i>0$. Since $s\ge\Delta+6$,
the generic term is bounded by
\[
 16q^{\Delta-s+1}\le16\,2^{-5}=\frac12.
\]
If $\mathcal M>0$, the choice of $L_s$ ensures
\[
 s(B_{d,c}-n\varepsilon_{d,c}\ell)\le-\log_2(64\mathcal M),
\]
so each dense contribution is at most $1/(4\mathcal M)$, uniformly for
$q\ge2$. The elementary estimates used for $G_{q,s,\ell}$ give
$G^{(i)}_{q,s,\ell}<16q^{\Delta_i+1-s-ns\varepsilon_{j_i,a_i}\ell}\le1/(4\mathcal M)$
for $1\le i\le\kappa$. When $\lambda>1$, we also have
\[
 R_{q,s,\ell}<q^{-ns(\lambda-1)\ell}\le\frac1{4\mathcal M}.
\]
There are exactly $\mathcal M$ terms of these three types, counted with
repetition. The finite failure bound is therefore at most $3/4$.
If $\mathcal M=0$, only the generic contribution remains and it is at most
$1/2$. In either case the required strict inequality below one
holds for every $\ell$ under consideration.
\end{proof}

\begin{remark}\label{rem:unique-ratio}
The strict inequality between the ratios in
Theorem~\ref{thm:simultaneous} is part of the result. If an additional
ratio equals $\rho_0$, some coefficients $\alpha c-d$ may vanish,
and our bound for the additional failure probability then no
longer tends to zero.
Moreover, an arbitrary simultaneous list obtained here is asserted
only in ambient ranks $r=s\ell$. The next section obtains all large
admissible ranks by selecting a particular profile that is preserved
under direct sums.
\end{remark}

The geometric construction always gives a nondegenerate code $C_U$.
For later applications we also need its dual to be nondegenerate.
This is a separate requirement, controlled by one additional event.

\begin{lemma}
\label{lem:kernel-nondegenerate}
Let $A\in\mathbb F_{q^{ns}}^{\ell\times T}$, and let
$\beta=(\beta_1,\ldots,\beta_s)$ be an $\mathbb F_q$-basis of
$\mathbb F_{q^s}$, also viewed as an $\mathbb F_{q^n}$-basis of
$\mathbb F_{q^{ns}}$. Let
$G\in\mathbb F_{q^n}^{s\ell\times sT}$ represent the
$\mathbb F_{q^n}$-linear map
$A:\mathbb F_{q^{ns}}^T\to\mathbb F_{q^{ns}}^\ell$
in the domain basis $(\beta_b\mathbf e_c)_{c,b}$ and codomain
basis $(\beta_b\mathbf f_u)_{u,b}$, where $(\mathbf e_c)$ and
$(\mathbf f_u)$ are the respective standard coordinate bases. Put
\[
 D=\ker_{\mathbb F_{q^n}}G,
 \qquad L=\operatorname{row}_{\mathbb F_{q^{ns}}}(A).
\]
Then $D$ is nondegenerate if and only if
$L\cap\mathbb F_{q^s}^{T}=0$.
\end{lemma}

\begin{proof}
Let $\beta^*=(\beta_1^*,\ldots,\beta_s^*)$ be the basis dual to
$\beta$ for the trace from $\mathbb F_{q^{ns}}$ to
$\mathbb F_{q^n}$. Its elements belong to $\mathbb F_{q^s}$, because
the restriction of this trace to $\mathbb F_{q^s}$ equals the trace
from $\mathbb F_{q^s}$ to $\mathbb F_q$. Write $\iota_\beta$ and
$\iota_{\beta^*}$ for coordinate expansion on
$\mathbb F_{q^{ns}}^T$ in these two bases. For $w,v$ in that space,
\[
 \operatorname{Tr}_{\mathbb F_{q^{ns}}/\mathbb F_{q^n}}(w\cdot v)
 =\iota_{\beta^*}(w)\cdot\iota_\beta(v).
\]
Since $\ker_{\mathbb F_{q^{ns}}}A$ is
$\mathbb F_{q^{ns}}$-linear, its trace annihilator equals its
ordinary dot-product annihilator: apply the trace identity to all
scalar multiples of a vector and use nondegeneracy of the trace.
The ordinary annihilator is $L$, so
\[
 D^\perp=\iota_{\beta^*}(L).
\]
A nonzero vector in this expanded row space has rank one over
$\mathbb F_q$ precisely when it is the expansion of
$ez$, for some $e\in\mathbb F_{q^n}^*$ and
$0\ne z\in\mathbb F_{q^s}^T$. Indeed, all of its coordinates must
be $\mathbb F_q$-multiples of $e$, and $\beta^*$ is an
$\mathbb F_q$-basis of $\mathbb F_{q^s}$.
As $L$ is $\mathbb F_{q^{ns}}$-linear, such a vector exists if and
only if $L\cap\mathbb F_{q^s}^T\ne0$. The dual criterion for
nondegeneracy from Section~\ref{sec:preliminaries} completes the proof.
\end{proof}

\begin{proposition}
\label{prop:nondegenerate-kernel}
Under the assumptions of Theorem~\ref{thm:simultaneous}, suppose
additionally that $(n-1)j>a$, and define
\[
 \mathcal N_{q,s,\ell}
 =\gauss{T}{1}{q^s}
   \frac{q^{ns\ell}-1}{q^{nsT}-1}.
\]
If the finite failure bound in that theorem, with $\mathcal N_{q,s,\ell}$
added, is less than one, the system $U$ can be chosen with
$C_U^\perp$ nondegenerate. If $16q^{nh_0j-a^2-s+1}<1$, this
additional requirement can be imposed for every sufficiently
large admissible $\ell$.

More explicitly, in Proposition~\ref{prop:uniform-threshold},
replace $\mathcal M$ by $\mathcal M+1$ and include, with this new
value of $\mathcal M$, the further integer
\[
 \left\lceil
 \frac{s^{-1}\log_2(16\mathcal M)-1}{\gamma}
 \right\rceil,
 \qquad \gamma=(n-1)\frac{nj}{a}-n>0,
\]
in the maximum defining $L_s$. The enlarged threshold ensures
nondegeneracy uniformly in $q$ as well as all the stated
intersection constraints. Nondegeneracy of the associated kernel
codes is preserved under direct sums of the systems.
\end{proposition}

\begin{proof}
Use the random matrix of Theorem~\ref{thm:simultaneous}.
Right multiplication by $\operatorname{GL}_T(\mathbb F_{q^{ns}})$
preserves the uniform matrix distribution and the full-row-rank
event, and acts transitively on the $\ell$-dimensional row spaces.
Conditional on full row rank, $L$ is therefore uniform among the
$\ell$-dimensional $\mathbb F_{q^{ns}}$-subspaces of
$\mathbb F_{q^{ns}}^T$. There are exactly
$\gauss{T}{1}{q^s}$ lines over $\mathbb F_{q^{ns}}$ containing a
nonzero vector of $\mathbb F_{q^s}^T$: two vectors in
$\mathbb F_{q^s}^T$ generate the same such line if and only if
they are proportional over $\mathbb F_{q^s}$. A fixed line is
contained in $L$ with conditional probability
\[
 \frac{q^{ns\ell}-1}{q^{nsT}-1}.
\]
Lemma~\ref{lem:kernel-nondegenerate} and a union bound show that
the probability of full row rank together with a degenerate
expanded kernel is at most $\mathcal N_{q,s,\ell}$. On the existing good
events, the columns of the expanded matrix form an
$\mathbb F_q$-basis of $U$, so that kernel is $C_U^\perp$.
Adding this event to the previous failure bound is therefore
valid; no independence between the events is required.

The two elementary estimates
\[
 \gauss{T}{1}{q^s}<2q^{s(T-1)},\qquad
 \frac{q^{ns\ell}-1}{q^{nsT}-1}<2q^{ns(\ell-T)}
\]
give
\[
 \mathcal N_{q,s,\ell}
 <4q^{s[n\ell-(n-1)T-1]}
 =4q^{-s(\gamma\ell+1)}.
\]
Here $(n-1)j>a$ implies $nj>a$ and $\gamma>0$, so this term
tends to zero. Under the uniform auxiliary hypothesis
$s\ge\Delta+6$, the extra ceiling makes it at most $1/(4\mathcal M)$
for every $q\ge2$. The dense terms, the full-span terms of the
additional constraints and the row-rank term have
the same bounds after recomputing the other ceilings with the
new value of $\mathcal M$. Their total together with the generic term
is at most $3/4$.

Finally, choose the concatenation of bases for a direct sum of
systems. Its coordinate matrix is block diagonal, so its kernel
code is the direct sum of the individual kernel codes. Rank
supports in the separate coordinate blocks form a direct sum.
If each kernel has full rank support in its block, their sum has
full rank support in all coordinates. This proves the last claim.
\end{proof}

\section{Stable profiles and all sufficiently large admissible ranks}
\label{sec:all-r}

Fix again $1\le h\le k$ with $\gcd(h,k+1)=1$, and put
$a=k+1$, $j=a-h$, assuming $nj\ge a$.
Define a function on $\{0,\ldots,h\}$ by
\[
 f(0)=0,\qquad
 f(d)=\left\lfloor\frac{ad}{h}\right\rfloor\quad(1\le d<h),
 \qquad f(h)=a-1.
\]
We say that $U\le_{\mathbb F_q}V$, with $V$ an $\mathbb F_{q^n}$-space, satisfies
the \emph{profile $f$} if
\[
 \dim_{\mathbb F_q}(U\cap H)\le f(d)
 \qquad\text{whenever }H\le_{\mathbb F_{q^n}} V,\quad\dim_{\mathbb F_{q^n}} H=d\le h.
\]
In particular, the constraint at $d=h$ is $(h,k)$-evasivity.

\begin{lemma}\label{lem:profile-sum}
If $U_1\le_{\mathbb F_q}V_1$ and $U_2\le_{\mathbb F_q}V_2$ satisfy
the profile $f$, then $U_1\oplus U_2$ satisfies it in
$V_1\oplus V_2$. The same conclusion holds for any finite direct sum.
\end{lemma}

\begin{proof}
Let $H\le_{\mathbb F_{q^n}} V_1\oplus V_2$ have dimension $d\le h$, and let $\pi_1$
be projection onto $V_1$. Put
\[
 d_1=\dim_{\mathbb F_{q^n}}\pi_1(H),\qquad d_2=\dim_{\mathbb F_{q^n}}(H\cap V_2),
 \qquad d_1+d_2=d.
\]
The image under $\pi_1$ of $(U_1\oplus U_2)\cap H$ is contained in
$U_1\cap\pi_1(H)$, and its kernel is $U_2\cap H\cap V_2$.
The dimension formula over $\mathbb F_q$ gives
\[
 \dim_{\mathbb F_q}((U_1\oplus U_2)\cap H)\le f(d_1)+f(d_2).
\]
If $d<h$, the right side is at most
$\lfloor a(d_1+d_2)/h\rfloor=f(d)$, by the elementary floor
inequality. If $d=h$ and one of $d_1,d_2$ is zero, the conclusion
is immediate. If $d=h$ and both are positive, coprimality implies
$ad_1/h\notin\mathbb Z$, and hence
\[
 f(d_1)+f(d_2)
 =\left\lfloor\frac{ad_1}{h}\right\rfloor
  +\left\lfloor a-\frac{ad_1}{h}\right\rfloor
 =a-1=f(h).
\]
This also covers $h=1$, since then the last case cannot occur.
Induction proves the assertion for a finite number of summands.
\end{proof}

\begin{theorem}\label{thm:all-r}
Let
\[
 1\le h\le k,\qquad\gcd(h,k+1)=1,
 \qquad a=k+1,\quad j=a-h,\quad
 \nu=\frac{a}{\gcd(n,a)}.
\]
Assume $nj\ge a$. There is a computable integer $r_0=r_0(n,h,k)$ such that, for every
prime power $q$ and every $r\ge r_0$ with $\nu\mid r$, there is a
spanning $(h,k)$-evasive system $U\le \mathbb F_{q^n}^r$ with
\[
 \dim_{\mathbb F_q} U=\frac{rnj}{a}.
\]
The system can be chosen to satisfy the entire profile $f$ above.
If $k<n$, it is maximum. If $(n-1)j>a$, it may additionally be
chosen so that $C_U^\perp$ is nondegenerate, with a threshold still
independent of $q$.
The divisibility condition $\nu\mid r$ is equivalent to the
integrality of $rnj/a$.
\end{theorem}

\begin{proof}
First impose the profile on one auxiliary-field block. In addition
to the primary pair $(h,k)$, consider the pairs $(d,f(d))$ for
$1\le d<h$. These satisfy
\[
 d\le f(d)\le a-2.
\]
The lower inequality follows from $a>h$. For the upper inequality,
$ad/h\le a-a/h<a-1$. Put $a_d=f(d)+1$ and $j_d=a_d-d$. Then
$a_dh>ad$, which is equivalent to $a_d(a-h)<a(a_d-d)$. Since
$j_d=a_d-d$ and $j=a-h$ are positive, this gives
\[
 \frac{a_d}{j_d}<\frac{a}{a-h}=\frac aj.
\]
The primary pair is therefore the unique pair of largest ratio in
this finite list. Theorem~\ref{thm:simultaneous} and
Proposition~\ref{prop:uniform-threshold} construct systems
with profile $f$ in rank $s\ell$, for every sufficiently large
admissible $\ell$, with a threshold uniform in $q$.

We give an explicit selection procedure for two blocks. Put
$\Delta=nhj-a^2$, and choose distinct primes $s_1<s_2$ greater than
$\max\{n,\Delta+5\}$. They are coprime to $n$ and satisfy the
uniform auxiliary condition $s_i\ge\Delta+6$. For the profile list
just defined, compute the corresponding integers $L_{s_i}$ from
Proposition~\ref{prop:uniform-threshold}. If $(n-1)j>a$ and
nondegenerate kernel codes are required, use the enlarged thresholds
of Proposition~\ref{prop:nondegenerate-kernel} instead.
Choose a prime $p_1$
greater than
\[
 \max\{s_1,s_2,L_{s_1}/\nu\},
\]
and a prime $p_2$ greater than
\[
 \max\{p_1,s_2,L_{s_2}/\nu\}.
\]
For definiteness, the least prime satisfying each requirement may
be used, so every choice is computable by finite search. Set
\[
 \ell_1=\nu p_1,\qquad \ell_2=\nu p_2,
 \qquad u_1=s_1p_1,\qquad u_2=s_2p_2.
\]
Since $a\mid n\nu$ and $\ell_i\ge L_{s_i}$, the construction gives,
for every $q$, profile-$f$ systems in the two ranks
$r_1=\nu u_1$ and $r_2=\nu u_2$. The four primes are distinct, so
$\gcd(u_1,u_2)=1$.

Every integer $R\ge(u_1-1)(u_2-1)$ is a nonnegative integer combination
of $u_1$ and $u_2$. Indeed, choose $0\le x\le u_2-1$ with
$xu_1\equiv R\pmod {u_2}$. Then $R-xu_1\ge-(u_2-1)$; as this is a multiple
of $u_2$, it is nonnegative. Thus $R=xu_1+yu_2$ with $x,y\ge0$.
For $r=\nu R$, take the direct sum of $x$ copies of a system of
rank $r_1$ and $y$ copies of a system of rank $r_2$.
Lemma~\ref{lem:profile-sum} preserves the profile, and the summands
span their respective $\mathbb F_{q^n}$-spaces. The new system therefore spans
$\mathbb F_{q^n}^r$ and has dimension
\[
 x\frac{r_1nj}{a}+y\frac{r_2nj}{a}
 =\frac{rnj}{a}.
\]
When requested, nondegeneracy of $C_U^\perp$ is preserved by the
same direct sums, by Proposition~\ref{prop:nondegenerate-kernel}.
If $k<n$, Theorem~\ref{thm:evasive-bound} proves maximality.
One explicit permissible
threshold is
\[
 r_0=\max\{h+1,\ \nu(u_1-1)(u_2-1)\}.
\]
All choices preceding this formula depend only on $n,h,k$.
Finally, $\gcd(a,j)=1$ implies
$a\mid rnj$ if and only if $a\mid rn$, which is equivalent to
$\nu\mid r$.
\end{proof}

\begin{remark}\label{rem:scope-all-r}
For $h=k\le n-2$, the existence conclusion is already included in
the stronger explicit results of Giannoni and Neri~\cite{GN}. The
endpoint $h=k=n-1$ is elementary: $\mathbb F_q^r\le \mathbb F_{q^n}^r$ is a
spanning maximum $(n-1)$-scattered subspace. For $h=1<k$ the
conclusion is known or elementary. Take duals with respect to
$(x,y)\mapsto\Tr_{q^n/q}(x\cdot y)$. If $W\le\mathbb F_{q^n}^{\,r}$ has
$\mathbb F_q$-dimension $t$ and $P$ is a one-dimensional
$\mathbb F_{q^n}$-subspace, then
$\dimq(W^\perp\cap P)=\dimq(W\cap P^\perp)+n-t$. For $r=k+1$ and
$W=\{(x,x^q,\ldots,x^{q^k}):x\in\mathbb F_{q^n}\}$, which has dimension
$n$ and meets every hyperplane in the kernel of a nonzero $q$-polynomial
of $q$-degree at most $k$, the subspace $W^\perp$ is therefore
$(1,k)$-evasive of the maximum dimension $kn$. It spans, because
otherwise $W$ would contain a one-dimensional $\mathbb F_{q^n}$-subspace.
Direct sums, which preserve $(1,k)$-evasivity by
Lemma~\ref{lem:profile-sum}, give every $r$ with $(k+1)\mid r$; for
$k=n-2$ these are all admissible ranks. For $k=n-1$, the dual of
$\mathbb F_q^r$ works for every $r$. For $k\le n-3$, the dual of a
maximum $k$-scattered subspace is maximum $(1,k)$-evasive
\cite[Theorem~2.15]{GRSZ}, and it spans for the same reason, so
\cite{GN} gives every sufficiently large admissible rank. Hence the
theorem is new only for $2\le h<k$. The argument above
extends the auxiliary-field count to $h<k$ and uses a stable
profile to pass from auxiliary ranks to every sufficiently large
admissible rank. Both uses of $\gcd(h,k+1)=1$ are visible in the
proof: it excludes equality in the dense strata of the primary
count, and it gives the strict floor identity at the boundary of
the direct-sum profile.
\end{remark}

\section{Generalized weights and optimal minimal codes}\label{sec:coding}

\subsection{A threshold bound and its attainment}
The following bound is a direct consequence of the evasivity dictionary
of \cite{MNT} and Theorem~\ref{thm:evasive-bound}. It identifies the
extremal coding parameters reached by our construction. Generalized
rank weights have also been used to distinguish scattered-space
constructions and their associated codes \cite{ZZ,BGM}.

\begin{proposition}\label{prop:weight-bound}
Let $D\le \mathbb F_{q^n}^N$ have dimension $K$, and assume that $D^\perp$ is nondegenerate. If
\[
 1\le j\le K,\qquad j+2\le d_0\le n+1,\qquad d_j(D)\ge d_0,
\]
then
\begin{equation}\label{eq:weight-bound}
 jnK\le N(jn-d_0+1).
\end{equation}
\end{proposition}
\begin{proof}
Put $r=N-K$, $h=d_0-1-j$, and $k=d_0-2$. By \eqref{eq:generalized-singleton},
\[
 1\le h\le r-1,\qquad h\le k<n.
\]
Since $d_0\le d_j(D)\le N$, we also have $k+1\le N$. The nondegenerate code $D^\perp$ defines a spanning $q$-system $U\le \mathbb F_{q^n}^r$ of dimension $N$. Proposition~\ref{prop:dictionary} says that $U$ is $(h,k)$-evasive. Therefore
\[
 N\le\frac{rnj}{d_0-1}.
\]
Substituting $r=N-K$ and rearranging proves \eqref{eq:weight-bound}.
\end{proof}

\begin{corollary}\label{cor:weight-attainment}
Fix $n,j,d_0$ with
\[
 1\le j,\qquad j+2\le d_0\le n+1,\qquad
 \gcd(j,d_0-1)=1,\qquad nj>d_0-1.
\]
For every sufficiently large $r$ satisfying $(d_0-1)\mid rn$, and every prime power $q$, there is an $\mathbb F_{q^n}$-linear code $D$ with
\[
 N=\frac{rnj}{d_0-1},\qquad K=N-r,\qquad
 d_j(D)\ge d_0,
\]
whose dual is nondegenerate and which attains \eqref{eq:weight-bound}. The threshold for $r$ can be chosen independently of $q$. If $d_0\le n$, then $d_j(D)=d_0$.
\end{corollary}
\begin{proof}
Apply Theorem~\ref{thm:all-r} with
\[
 h=d_0-1-j,\qquad k=d_0-2.
\]
Its hypotheses hold because $\gcd(h,k+1)=\gcd(j,d_0-1)=1$. The maximum $q$-system $U$ provided by that theorem has dimension $N=rnj/(d_0-1)$. Put $D=C_U^\perp$. The strict inequality $nj>d_0-1$ ensures that $K=N-r$ grows with $r$, so $j\le K$ for sufficiently large $r$. Proposition~\ref{prop:dictionary} gives $d_j(D)\ge d_0$, and the parameter identity gives equality in \eqref{eq:weight-bound}.

If $d_0\le n$ and $d_j(D)\ge d_0+1$, Proposition~\ref{prop:weight-bound} applied at threshold $d_0+1$ would give
\[
 jnK\le N(jn-d_0),
\]
contradicting the equality already obtained.
\end{proof}

\begin{remark}\label{rem:endpoint}
At the endpoint $d_0=n+1$, the corollary asserts the threshold $d_j(D)\ge n+1$. It does not assert equality of the generalized weight with $n+1$: the next threshold lies outside the range of Proposition~\ref{prop:weight-bound}. The excluded case $nj=d_0-1$ is exactly $j=1$, $d_0=n+1$, where the displayed parameters give $K=0$ and no first generalized weight is defined.
\end{remark}

\subsection{Minimality and the universal length bound}
For the rest of this section, as in \cite{ABNR,BMN,XKH}, we write $m$
for the extension degree: codes are linear over $\mathbb F_{q^m}$, and
the results of Sections~\ref{sec:preliminaries}--\ref{sec:all-r} are
applied with $n=m$. For a nonzero codeword $c$, its rank support is unchanged by multiplication by an element of $\mathbb F_{q^m}^*$. An $\mathbb F_{q^m}$-linear code $C$ is \emph{minimal in the rank metric} if, for nonzero $c,c'\in C$,
\[
 \Supp_q(c')\subseteq\Supp_q(c)\quad\Longrightarrow\quad c'\in \mathbb F_{q^m}c.
\]
Minimal rank-metric codes and their geometric interpretation by linear cutting blocking sets were developed in \cite{ABNR}. The characterization needed here is that, for codes with $\dim_{\mathbb F_{q^m}} C\ge2$,
\begin{equation}\label{eq:minimal-weight}
 C\text{ minimal}\quad\Longleftrightarrow\quad d_2(C)\ge m+1.
\end{equation}
For nondegenerate codes this is \cite[Theorem 3.4]{BMN}; the formulation allowing degeneracy follows by deleting zero coordinates after a rank isometry, and is also included in the higher-minimality characterization of
\cite[Theorem~4.4]{XKH}.

Let $\varpi_{q^m/q}(K)$ denote the minimum length of an $\mathbb F_{q^m}$-linear minimal rank-metric code of dimension $K$. Allowing degenerate codes does not change this minimum, since deleting degenerate coordinates preserves minimality and dimension.

\begin{proposition}\label{prop:minimal-lower}
If $m\ge3$ and $C$ is an $[N,K]_{q^m/q}$ minimal code with $K\ge2$, then
\[
 N\ge\max\{2K,K+m-1\}.
\]
\end{proposition}
\begin{proof}
The bound $N\ge K+m-1$ is \cite[Corollary 5.10]{ABNR};
it also follows directly from
$m+1\le d_2(C)\le N-K+2$.
It remains to prove $N\ge2K$; for $m=3$ this is
\cite[Theorem~5.3(1)]{XKH}.
If $d_1(C)\ge2$, then $C^\perp$ is nondegenerate. Apply Proposition~\ref{prop:weight-bound} with $n=m$, $j=2$, $d_0=m+1$, using \eqref{eq:minimal-weight}. It gives $2mK\le Nm$.

Suppose instead that $d_1(C)=1$. A rank isometry, followed by scaling a generator, puts a rank-one word in the form $\mathbf e_1$. Eliminating the first coordinate of the remaining generators gives
\[
 C=\mathbb F_{q^m}\mathbf e_1\oplus(0,C'),\qquad \dim_{\mathbb F_{q^m}} C'=K-1.
\]
For every nonzero $v\in C'$, the two-dimensional subcode spanned by $\mathbf e_1$ and $(0,v)$ has rank-support weight $1+\rk_q(v)$. By \eqref{eq:minimal-weight}, $\rk_q(v)=m$. Thus $C'$ has minimum rank distance $m$ and is nonzero, so $N-1\ge m$. The matrix Singleton bound gives
\[
 m(K-1)\le N-1.
\]
To verify it, expand codewords into $m\times(N-1)$ matrices and delete any $m-1$ rows. This operation is injective on $C'$, since a nonzero difference supported on the deleted rows would have rank at most $m-1$. Only $N-1$ entries remain, proving the displayed inequality. Finally,
\[
 1+m(K-1)-2K=(m-2)K-m+1\ge m-3\ge0.
\]
\end{proof}

\noindent
The two terms agree at $K=m-1$: the first dominates for
$K\ge m-1$, and the second for $K\le m-1$.
The assumption $m\ge3$ is necessary for the bound $2K$,
since $\varpi_{q^2/q}(K)=2K-1$ \cite[Corollary~5.1(2)]{XKH}.
Further lower bounds for particular parameter ranges are given
in \cite[Theorem~5.2 and Corollaries~5.1 and~5.3]{XKH}. For instance,
\cite[Corollary~5.3(1)]{XKH} gives $N\ge K+m$ whenever
$K\ge\lfloor\sqrt m\rfloor+2$, which exceeds $\max\{2K,K+m-1\}$ for
$K<m$. In particular $\varpi_{q^m/q}(m-1)\ge2m-1$ for $m\ge5$, which
also follows directly: a minimal $[2m-2,m-1]_{q^m/q}$
code is nondegenerate by Proposition~\ref{prop:minimal-lower}, and by
\cite[Theorem~3.3]{BMN} its system would be an $(m-3)$-scattered
subspace of $\mathbb F_{q^m}^{\,m-1}$ of dimension $2m-2>(m-1)m/(m-2)$,
contrary to Theorem~\ref{thm:evasive-bound}. For $m\ge4$ we have not
found the bound $N\ge2K$ in the literature.

\begin{theorem}\label{thm:minimal-optimal}
For every odd integer $m\ge3$, there is an integer $K_0(m)$ such that
\[
 \varpi_{q^m/q}(K)=2K
\]
for every prime power $q$ and every $K\ge K_0(m)$.
\end{theorem}
\begin{proof}
Apply Theorem~\ref{thm:all-r} with $n=m$, $h=m-2$, and $k=m-1$. The coprimality condition is
\[
 \gcd(m-2,m)=\gcd(2,m)=1,
\]
and the divisibility condition $m\mid rm$ imposes no restriction on $r$. Hence, for every sufficiently large $K$, uniformly in $q$, there is a spanning maximum $(m-2,m-1)$-evasive subspace
\[
 U\le\mathbb F_{q^m}^{\,K},\qquad \dimq U=2K.
\]
The code $D=C_U^\perp$ has length $2K$ and dimension $K$. Proposition~\ref{prop:dictionary}, with $j=2$, gives $d_2(D)\ge m+1$, so $D$ is minimal by \eqref{eq:minimal-weight}. Proposition~\ref{prop:minimal-lower} proves optimality.

The constructed code is also nondegenerate. Otherwise, deleting degenerate coordinates after a rank isometry would yield a minimal code of dimension $K$ and length less than $2K$, contrary to Proposition~\ref{prop:minimal-lower}.
\end{proof}

For $m=3$, the equality $\varpi_{q^3/q}(K)=2K$ is already
known for every $K\ge2$ by \cite[Corollary~5.4(2)]{XKH}. For arbitrary
$m\ge3$, that paper gives the upper bound $2K+m-3$
\cite[Corollary~5.4(1)]{XKH}.
Theorem~\ref{thm:minimal-optimal} therefore determines the eventual
value in the minimum-length problem of
\cite[Remark 6.12]{ABNR} for every fixed odd $m\ge5$.
This parameter range is not covered by the isolated
short-code constructions of \cite{BMN,LLMT}.

The theorem also determines every possible length in this range.
Let $m$ be odd and $K\ge K_0(m)$. Proposition~\ref{prop:minimal-lower}
excludes all lengths below $2K$. Conversely, append to a generator
matrix of a code from Theorem~\ref{thm:minimal-optimal} further columns
that are $\mathbb F_q$-independent of the previous ones. Projecting a
subcode onto the original coordinates maps its rank support onto the
rank support of a subcode of the same dimension, so the generalized
weights do not decrease and $d_2\ge m+1$ is preserved. Hence
nondegenerate minimal $[N,K]_{q^m/q}$ codes exist for every
$2K\le N\le Km$. In this range, this settles the existence question
for the lengths between the bounds of \cite{ABNR}, which is listed
among the open problems in \cite[Section~5.2]{Scotti}.

\subsection{The first two generalized weights}
The stable profile also determines the minimum rank distance of the
ordinary minimal codes above. With a choice of at least two summands,
it determines their second generalized weight as well.

\begin{corollary}\label{cor:minimal-exact-weights}
For each fixed odd $m\ge3$, the codes in
Theorem~\ref{thm:minimal-optimal} can be chosen, for every sufficiently
large $K$ uniformly in $q$, with
\[
 d_1=\frac{m+1}{2},\qquad d_2=m+1.
\]
The first value is the largest possible minimum rank distance of an
$\mathbb F_{q^m}$-linear $[2K,K]$ code for these $K$.
\end{corollary}
\begin{proof}
Write $m=2b+1$. Let $U\le\mathbb F_{q^m}^{\,K}$ have
$\mathbb F_q$-dimension $2K$ and satisfy the stable profile for
$(h,k)=(m-2,m-1)$. Put $C=C_U^\perp$.
If $b\ge2$, the profile at dimension $b-1$ gives
\[
 \dimq(U\cap H)\le
 \left\lfloor\frac{(2b+1)(b-1)}{2b-1}\right\rfloor=b-1
 \qquad\left(\dim_{\mathbb F_{q^m}}H=b-1\right).
\]
Thus $U$ is $(b-1)$-scattered, and
Proposition~\ref{prop:dictionary}, with $j=1$, gives $d_1(C)\ge b+1$.
For $b=1$, the inequality $d_1(C)\ge2$ follows from the
$\mathbb F_q$-independence of the columns defining $U$.
Conversely, expanding any $\mathbb F_{q^m}$-linear $[2K,K]$ code
into $m\times2K$ matrices and deleting $d_1-1$ rows gives
\[
 mK\le2K(m-d_1+1),
\]
so $d_1\le\lfloor(m+2)/2\rfloor=b+1$.

For the second weight, use the two fixed seed systems in the proof
of Theorem~\ref{thm:all-r}. For all sufficiently large $K$, their
direct-sum construction uses at least two summands: it suffices
to take $K$ larger than both seed ranks as well as the
semigroup threshold. The associated code is
$C=C_1\oplus\cdots\oplus C_z$, with $z\ge2$, on disjoint
coordinate blocks. The preceding argument applies to each seed
code, so choose $c_i\in C_i$ of rank $b+1$ for $i=1,2$.
The two words supported on the first and second blocks generate
a two-dimensional $\mathbb F_{q^m}$-subcode whose rank support
has dimension $2(b+1)=m+1$. Hence $d_2(C)\le m+1$.
The reverse inequality follows from minimality, proving equality.
\end{proof}

Viewed as $m\times2K$ matrices, the codes attain the rounded-down
matrix Singleton bound, and are therefore quasi-MRD in the sense of
\cite[Definition~10]{DGLR} (see also \cite{MNT}).
Since $m$ is odd, the unrounded matrix Singleton inequality is strict.
Scattered systems that are maximal under inclusion but not maximum,
together with the associated quasi-MRD and non-extendable MRD codes,
are studied in \cite{Gia26,BGMN26}.

\begin{example}\label{ex:weights}
For $(n,h,k)=(5,2,2)$, every sufficiently large $r$ divisible by $3$ yields a code with
\[
 (N,K)=\left(\frac{5r}{3},\frac{2r}{3}\right),\qquad d_1=4.
\]
This first family is MRD: it meets the matrix Singleton bound,
and its scattered construction is already covered by \cite{GN}.
For $(n,h,k)=(5,3,4)$, every sufficiently large $r$ yields a
minimal $[2r,r]_{q^5/q}$ code which can be chosen with
$d_1=3$ and $d_2=6$. Here equality at the endpoint of
Remark~\ref{rem:endpoint} follows from the direct-sum construction.
\end{example}

\subsection{The stable profile as a coding property}
For the ordinary minimal codes, the stable profile of
Section~\ref{sec:all-r} has a direct coding meaning. This explains its
stability under direct sums, and it shows that no weaker condition on
the blocks can be used in a construction by direct sums.

\begin{proposition}\label{prop:profile-coding}
Let $m=2b+1\ge5$, let $K\ge m-1$, and let $U\le\mathbb F_{q^m}^{\,K}$
be a spanning $q$-system with $\dimq U=2K$. Put $C=C_U^\perp$, a
$[2K,K]_{q^m/q}$ code, and let $f$ be the stable profile for
$(h,k)=(m-2,m-1)$.
\begin{enumerate}
\item One has $f(d)=d$ for $0\le d\le b-1$, and $f(d)=d+1$ for
$b\le d\le m-2$.
\item The system $U$ satisfies $f$ if and only if $C$ is minimal and
$d_1(C)=b+1$.
\item If $C_1$ and $C_2$ are minimal $[2K_i,K_i]_{q^m/q}$ codes with
$d_1(C_1)=d_1(C_2)=b+1$, then $C_1\oplus C_2$ is minimal with
$d_1=b+1$.
\item Conversely, if $C_1\oplus C_2$ is minimal, where each $C_i$ is a
nonzero $[2K_i,K_i]_{q^m/q}$ code, then $d_1(C_1)=d_1(C_2)=b+1$.
\end{enumerate}
\end{proposition}

\begin{proof}
For $1\le d<m-2=2b-1$ we have $f(d)=d+\lfloor 2d/(2b-1)\rfloor$.
Here $2d<2b-1$ exactly when $d\le b-1$, while $2b-1<2d<2(2b-1)$ when
$b\le d\le 2b-2$. Together with $f(m-2)=m-1$, this proves (1).

By (1), the constraint of $f$ at $d=b-1$ says that $U$ is
$(b-1)$-scattered, and the constraint at $d=m-2$ is
$(m-2,m-1)$-evasivity. Proposition~\ref{prop:dictionary}, applied with
$j=1$ and with $j=2$, translates them into $d_1(C)\ge b+1$ and
$d_2(C)\ge m+1$; its hypotheses hold because $b\ge2$ and $K\ge m-1$.
By \eqref{eq:minimal-weight}, the second inequality is minimality.
Every $[2K,K]_{q^m/q}$ code with $2K\ge m$ has $d_1\le b+1$, by the
matrix Singleton bound as in the proof of
Corollary~\ref{cor:minimal-exact-weights}. Hence the profile implies
the two coding conditions. Conversely, if $C$ is minimal and
$d_1(C)=b+1$, then $U$ is $(b-1)$-scattered and $(m-2,m-1)$-evasive.
Suppose that an $\mathbb F_{q^m}$-subspace $H$ of dimension $d$
violated the profile. Add vectors of $U$ outside $H$ one at a time,
which is possible because $U$ spans; each step raises both
$\dim_{\mathbb F_{q^m}}H$ and $\dimq(U\cap H)$ by at least one. Stopping
at dimension $b-1$ when $d\le b-1$, and at dimension $m-2$ when
$d\ge b$, contradicts one of the two conditions. This proves (2).

For (3), no $C_i$ contains a word of rank one, so $C_i^\perp$ is
nondegenerate and $C_i=C_{U_i}^\perp$ for a spanning $q$-system
$U_i\le\mathbb F_{q^m}^{\,K_i}$ of dimension $2K_i$. Moreover
$K_i\ge m-1$: a minimal code of dimension $K_i\ge2$ has length at least
$K_i+m-1$ by Proposition~\ref{prop:minimal-lower}, and a
one-dimensional $[2,1]_{q^m/q}$ code has $d_1\le2<b+1$. By (2), both
$U_i$ satisfy $f$. By Lemma~\ref{lem:profile-sum}, so does
$U_1\oplus U_2$, whose kernel code is $C_1\oplus C_2$, and (2) applies
again.

For (4), choose $c_i\in C_i$ of rank $d_1(C_i)$. The words $(c_1,0)$ and
$(0,c_2)$ span a two-dimensional subcode whose rank support has
dimension $d_1(C_1)+d_1(C_2)$. Minimality gives
$d_1(C_1)+d_1(C_2)\ge m+1=2b+2$. The matrix Singleton bound gives
$d_1(C_i)\le b+1$ for every $[2K_i,K_i]_{q^m/q}$ code, so both minimum
distances are equal to $b+1$.
\end{proof}

Thus Lemma~\ref{lem:profile-sum}, in this case, states that minimal
$[2K,K]_{q^m/q}$ codes with the largest possible minimum distance are
closed under direct sums, and item (4) shows that a direct-sum
construction of minimal codes of length $2K$ must use blocks of this
kind. Indeed, let $C_1\oplus\cdots\oplus C_z$, with $z\ge2$ nonzero
blocks, be a minimal code of length $2K$ and dimension $K$. Every block
is minimal, because its words are the words of the sum supported on
its coordinates. A block of dimension $K_i\ge2$ therefore has length at
least $2K_i$ by Proposition~\ref{prop:minimal-lower}, and a block of
dimension one has length at least two unless its length is one. A block
of length one would force $d_1(C_j)\ge m$ for every other block, hence
length at least $mK_j$ by the matrix Singleton bound, and then
$2K\ge1+m(K-1)$, which is impossible for $m\ge4$ and $K\ge2$. Since the
lengths add up to $2K$, every block has length twice its dimension, and
(4) applies to any two of them. Any such code, for a single further dimension and for some or all
$q$, can therefore be added to the seeds of
Corollary~\ref{cor:minimal-degree-five} below.

The minimality part of (3) also follows from \cite[Proposition~3.8]{BMN}.
In the language of codes, that result shows that $C_1\oplus C_2$ is
minimal whenever $C_1$ and $C_2$ are nondegenerate minimal codes with
$d_1(C_1)+d_1(C_2)\ge m+1$; the proof of (4) shows that this inequality
is also necessary.

For even $m$ the same reasoning identifies the relevant class, and it
shows where the coprimality obstruction enters.

\begin{proposition}\label{prop:even-degree}
Let $m=2b\ge4$, and define
\[
 g(d)=d\quad(0\le d\le b-1),\qquad g(d)=d+1\quad(b\le d\le m-2).
\]
\begin{enumerate}
\item If $C_1\oplus\cdots\oplus C_z$ is minimal, with $z\ge2$ and each
$C_i$ a nonzero $[2K_i,K_i]_{q^m/q}$ code, then $d_1(C_i)\le b+1$ for
every $i$, with equality for all indices except at most one. Equality
holds exactly when $2K_i\ge m$ and $C_i$ is MRD.
\item Let $K\ge m-1$, and let $U\le\mathbb F_{q^m}^{\,K}$ be a spanning
$q$-system with $\dimq U=2K$. Then $U$ satisfies the profile $g$, in
the sense of Section~\ref{sec:all-r}, if and only if $C_U^\perp$ is
minimal and MRD.
\item The profile $g$ is preserved by direct sums. Consequently, direct
sums of minimal MRD $[2K_i,K_i]_{q^m/q}$ codes with $K_i\ge2$ are
minimal MRD codes.
\end{enumerate}
In the language of Theorem~\ref{thm:simultaneous}, the MRD condition
in (2) is the constraint $(b-1,b-1)$, whose ratio $b$ equals the ratio
$m/2$ of the primary constraint $(m-2,m-1)$; moreover
$\gcd(m-2,m)=2$, so the primary constraint is not coprime either.
\end{proposition}

\begin{proof}
For (1), when $2K_i\ge m$ the matrix Singleton bound gives
$mK_i\le2K_i(m-d_1(C_i)+1)$, hence $d_1(C_i)\le b+1$, with equality
exactly in the MRD case. When $2K_i<m$ it gives
$d_1(C_i)\le K_i+1\le b$. As in
Proposition~\ref{prop:profile-coding}(4), minimum-rank words in two
different blocks give $d_1(C_i)+d_1(C_{i'})\ge m+1=2b+1$ for $i\ne i'$.
Hence at most one index has $d_1(C_i)\le b$.

For (2), the constraint of $g$ at $d=b-1$ is $(b-1)$-scattering, which
is equivalent to $d_1(C_U^\perp)\ge b+1$ by
Proposition~\ref{prop:dictionary}, and hence to the MRD property by the
Singleton bound, since $2K\ge m$. The constraint at $d=m-2$ is
$(m-2,m-1)$-evasivity, which is equivalent to minimality. The remaining
values of $g$ follow from these two constraints by adding vectors of
$U$ one at a time, as in the proof of
Proposition~\ref{prop:profile-coding}(2).

For (3), let $H\le V_1\oplus V_2$ have dimension $d\le m-2$, and let
$e_1,e_2$ be the dimensions called $d_1,d_2$ in the proof of
Lemma~\ref{lem:profile-sum}, so that
$e_1+e_2=d$ and $\dimq((U_1\oplus U_2)\cap H)\le g(e_1)+g(e_2)$. If
$e_1,e_2\le b-1$, the right side is $d\le g(d)$. If exactly one of them
is at least $b$, it is $d+1=g(d)$, because $d\ge b$. They cannot both
be at least $b$, since $d\le2b-2$. The statement on codes follows from
(2) as in Proposition~\ref{prop:profile-coding}(3): each block has no
word of rank one, and $K_i\ge m-1$ by
Proposition~\ref{prop:minimal-lower}, since $K_i\ge2$. Minimality of
these direct sums also follows from \cite[Proposition~3.8]{BMN}, as
noted after Proposition~\ref{prop:profile-coding}, because MRD blocks
have $d_1=b+1$.
\end{proof}

\begin{corollary}\label{cor:quartic}
Let $q$ be an odd power of $2$. For every integer $w\ge1$ there is a
minimal MRD $[8w,4w]_{q^4/q}$ code. In particular,
$\varpi_{q^4/q}(4w)=8w$.
\end{corollary}

\begin{proof}
For such $q$, Bartoli, Marino and Neri construct an $8$-dimensional
$\mathbb F_q$-subspace $U_0\le\mathbb F_{q^4}^{\,4}$ that is maximum
scattered and a linear cutting blocking set
\cite[Theorems~4.1 and~4.5]{BMN}. A
cutting blocking set spans the ambient space, and by
\cite[Theorem~3.3]{BMN} being cutting means being $(2,3)$-evasive.
Thus $U_0$ satisfies the profile $g$ of
Proposition~\ref{prop:even-degree} for $m=4$, namely $g(1)=1$ and
$g(2)=3$, and so does the direct sum $U_0^{\oplus w}$. By
Proposition~\ref{prop:even-degree}(2), its kernel code is a minimal MRD
$[8w,4w]_{q^4/q}$ code, and Proposition~\ref{prop:minimal-lower} shows
that its length is optimal. For $w=1$ the corollary is
\cite[Proposition~4.9 and Corollary~3.5]{BMN}, where the minimal MRD
code is $C_{U_0}$ itself rather than its dual; the new ingredient for
$w\ge2$ is the lower bound of Proposition~\ref{prop:minimal-lower}.
\end{proof}

\subsection{An explicit threshold in extension degree five}
\begin{corollary}
\label{cor:minimal-degree-five}
Let $\mathcal S_1=\{14,16,18,21,22,24,26,27,33\}$. The additive
semigroup generated by $\mathcal S_1$ consists of
\[
 14,\ 16,\ 18,\ 21,\ 22,\ 24,\ 26,\ 27,\ 28,\ 30
\]
and of all integers $K\ge32$. For every prime power $q$ and every $K$ in
this semigroup, there is a nondegenerate minimal $[2K,K]_{q^5/q}$ code
with $d_1=3$, and
\[
 \varpi_{q^5/q}(K)=2K .
\]
In particular, this holds for every $K\ge32$. If moreover
$K\in\{28,30,32\}$ or $K\ge34$, the code can be chosen with $d_2=6$.
\end{corollary}

\begin{proof}
Take $(n,h,k)=(5,3,4)$, so $a=5$, $j=2$, $\alpha=2/5$ and
$\Delta=5$. The stable profile is
\[
 f(0)=0,\qquad f(1)=1,\qquad f(2)=3,\qquad f(3)=4.
\]
Its additional constraints $(1,1)$ and $(2,3)$ both have ratio $2$,
smaller than $5/2$. In Theorem~\ref{thm:simultaneous}, the dense pair
$(1,2)$ occurs twice, from the primary constraint and from $(2,3)$, and
each additional constraint contributes a full-span term. The exponents
of $q$ in the bounds \eqref{eq:refined-counts} and
\eqref{eq:secondary-generic} are as follows:
\[
 \begin{array}{c|c|c}
 \text{term}&\text{multiplicity}&\text{exponent of }q\\ \hline
 \text{generic}&1&6-s\\
 D_{q,s,\ell}(1,2)&2&s(1-\ell)\\
 \text{full span of }(1,1)&1&2-s-s\ell\\
 \text{full span of }(2,3)&1&5-s-2s\ell
 \end{array}
\]
We construct a seed for every $s\ge7$ with $\gcd(s,5)=1$ and every
$\ell\ge2$. Here $T=2\ell\ge4$, so $sT\ge5$, $T\ge2$ for the dense pair,
and $s\ell\ge4$; for $\ell=2$ we use that
Theorems~\ref{thm:finite-evasive} and~\ref{thm:simultaneous} only
require $sT\ge a_i$. We also use $R_{q,s,\ell}<q^{-5s\ell}$. The
resulting bound is nonincreasing in $q\ge2$, $s\ge7$ and $\ell\ge2$:
every exponent is negative and decreases with $s$ and $\ell$, and every
product correction decreases as its argument increases. At
$(q,s,\ell)=(2,7,2)$ this bound is
\[
 \begin{aligned}
 \mathfrak B_1={}&
 \frac{C_2(2^5)}{2(1-2^{-7})}\prod_{u=2}^5(1-2^{-u})^{-1}
 +2^{-6}C_2(2^7)C_1(2^{35})\\
 &+\frac{2^{-19}(1-2^{-2})^{-1}C_1(2^5)
     +2^{-30}\prod_{u=2}^{4}(1-2^{-u})^{-1}C_2(2^5)}{1-2^{-7}}
 +2^{-70}<\frac{89}{100}.
 \end{aligned}
\]
For an explicit rational check, the generic and dense terms are
smaller than $8737/10000$ and $1575/100000$, and the remaining terms
together are smaller than $1/100000$; the sum of these bounds is
$88946/100000<89/100$. Thus every indicated seed has the full profile,
ambient rank $s\ell$ and $\mathbb F_q$-dimension $2s\ell$.

The elements of $\mathcal S_1$ are seed ranks: they are the products
$s\ell$ for the pairs $(s,\ell)$ equal to $(7,2)$, $(8,2)$, $(9,2)$,
$(7,3)$, $(11,2)$, $(8,3)$, $(13,2)$, $(9,3)$ and $(11,3)$.
Below $32$, the semigroup consists of the generators up to $27$ and of
$28=14+14$ and $30=14+16$. The following table certifies every rank from
$32$ through $45$:
\[
 \begin{array}{c|c@{\qquad}c|c@{\qquad}c|c}
 K&\text{decomposition}&K&\text{decomposition}&K&\text{decomposition}\\ \hline
 32&16+16&37&16+21&42&21+21\\
 33&33&38&16+22&43&21+22\\
 34&16+18&39&18+21&44&22+22\\
 35&14+21&40&18+22&45&21+24\\
 36&18+18&41&14+27&&
 \end{array}
\]
Adding copies of the generator $14$ covers every $K\ge32$. The integer
$31$ is not in the semigroup, since $31-g$ is not in it for any
generator $g\le31$; it is not a seed rank either, being prime. By
Lemma~\ref{lem:profile-sum}, matching these direct sums of seed systems
gives the full profile in every rank $K$ of the semigroup.

The kernel code $D=C_U^\perp$ has parameters $[2K,K]_{q^5/q}$ and
$d_2(D)\ge6$, by Proposition~\ref{prop:dictionary}.
Equation~\eqref{eq:minimal-weight} gives minimality, and
Proposition~\ref{prop:minimal-lower} gives optimality of its length.
It is nondegenerate, since deleting degenerate coordinates would
produce a shorter minimal code of the same dimension. The profile at
dimension one gives $d_1\ge3$, and the matrix Singleton bound gives
$5K\le2K(6-d_1)$, hence $d_1=3$.
Finally, every $K\ge28$ of the semigroup other than $33$ is a sum of at
least two generators: this holds for $28$, $30$, $32$ and the table
entries from $34$ to $45$, and for $K\ge46$ because $K-14\ge32$. Each
seed kernel has $d_1=3$, and a minimum-rank word in each of two
disjoint blocks produces a two-dimensional subcode of rank-support
dimension $6$. Hence $d_2\le6$, and $d_2=6$ for these $K$.
\end{proof}

\begin{remark}\label{rem:thresholds-odd}
The computation in Corollary~\ref{cor:minimal-degree-five} applies to
every odd $m=2b+1\ge5$, with $(n,h,k)=(m,m-2,m-1)$. Here
$\Delta=m^2-4m$, so the generic bound is below one only for
$s\ge m^2-4m+2$; this integer is always coprime to $m$. By
Proposition~\ref{prop:profile-coding}(1), every additional constraint
has $j_i\le2$, so the dense pairs are of the form $(1,c)$ with $c\le b$.
Their exponents increase with $c$, and the largest one is
\[
 E_\ell(1,b)=b^2-b-1-\ell,
\]
which is negative exactly when $\ell\ge b(b-1)=(m-1)(m-3)/4$. The
full-span terms have negative exponents for every $\ell\ge1$. For
$s\ge m^2-4m+2$ and $\ell\ge b(b-1)$, the resulting bound is therefore
nonincreasing in $q$, $s$ and $\ell$, as in
Corollary~\ref{cor:minimal-degree-five}, and it tends to zero as $s$ and
$\ell$ grow. Hence it is below one at some $(2,s_0,\ell_0)$, and every
pair $(s,\ell)$ with $s\ge s_0$, $\gcd(s,m)=1$ and $\ell\ge\ell_0$
supplies a seed of rank $s\ell$ for every $q$. The threshold $K_0(m)$
in Theorem~\ref{thm:minimal-optimal} may be taken to be the conductor of
the semigroup generated by these ranks; for $m=5$, with
$(s_0,\ell_0)=(7,2)$, this is Corollary~\ref{cor:minimal-degree-five}.
For the weights in Corollary~\ref{cor:minimal-exact-weights}, $K$ must
moreover be a sum of at least two seed ranks.
\end{remark}

\subsection{Higher minimality}
Let $1\le \tau<K$. An $[N,K]_{q^m/q}$ code $C$ is \emph{$\tau$-minimal}
if, for any two $\tau$-dimensional $\mathbb F_{q^m}$-subcodes $D,D'\le C$,
\[
 \Supp_q(D')\subseteq\Supp_q(D)
 \quad\Longrightarrow\quad \Supp_q(D')=\Supp_q(D).
\]
This is the definition of Xu, Kan and Han~\cite{XKH}, who write $r$ for
$\tau$; for $\tau=1$ it recovers ordinary minimality. Their
characterization \cite[Theorem~4.4]{XKH} is
\begin{equation}\label{eq:higher-minimal-weight}
 C\text{ is $\tau$-minimal}
 \quad\Longleftrightarrow\quad d_{\tau+1}(C)\ge m\tau+1.
\end{equation}
We write $\varpi_{q^m/q}(K,\tau)$ for the minimum length of a
$\tau$-minimal code of dimension $K$ over $\mathbb F_{q^m}$.
As for ordinary minimality, deleting degenerate coordinates preserves
the property and the dimension.

\begin{theorem}\label{thm:higher-minimal}
Fix integers $m\ge3$ and $1\le \tau\le m-2$ such that
$\gcd(m,\tau+1)=1$. There is an integer $K_0(m,\tau)$, independent of
$q$, such that, for every prime power $q$ and every $K\ge K_0(m,\tau)$,
there is a nondegenerate $\tau$-minimal code with parameters
\[
 [(\tau+1)K,K]_{q^m/q}.
\]
In particular, $\varpi_{q^m/q}(K,\tau)\le(\tau+1)K$ in this range.
\end{theorem}
\begin{proof}
Apply Theorem~\ref{thm:all-r} with
\[
 n=m,\qquad h=\tau(m-1)-1,\qquad k=m\tau-1.
\]
Then $a=m\tau$, $j=\tau+1$, and
\[
 \gcd(a,j)=\gcd(m,\tau+1)=1,\qquad nj>a,
 \qquad \frac{a}{\gcd(n,a)}=\tau.
\]
Moreover,
\[
 (n-1)j-a=(m-1)(\tau+1)-m\tau=m-\tau-1>0.
\]
By Proposition~\ref{prop:nondegenerate-kernel}, the systems in
Theorem~\ref{thm:all-r} can therefore be chosen with nondegenerate
associated kernel codes. Take $r=\tau K$. For every sufficiently large
$K$, uniformly in $q$, we obtain a spanning system
\[
 U\le\mathbb F_{q^m}^{\,\tau K},\qquad
 \dimq U=\frac{\tau K\,m(\tau+1)}{m\tau}=(\tau+1)K,
\]
which is $(\tau(m-1)-1,m\tau-1)$-evasive. The code $C=C_U^\perp$ has
length $(\tau+1)K$ and dimension $K$. Increase the threshold so that
$K\ge \tau+1$, $\tau K\ge h+1$, and $(\tau+1)K\ge m\tau$.
Proposition~\ref{prop:dictionary} now gives
$d_{\tau+1}(C)\ge m\tau+1$, and \eqref{eq:higher-minimal-weight} proves
$\tau$-minimality.
\end{proof}

The parameter $k=m\tau-1$ explains the need for the range $k\ge n$
in the geometric existence theorem. The dimension bound of
Theorem~\ref{thm:evasive-bound} is not used in this range.
The general existence bound of \cite[Theorem~5.4]{XKH}, for
$K\ge \tau+1$, is
\[
 \varpi_{q^m/q}(K,\tau)\le(\tau+1)K+\tau(m-\tau-2).
\]
Thus Theorem~\ref{thm:higher-minimal} removes the positive additive
term for $1\le \tau\le m-3$, subject to coprimality and for $K$
sufficiently large. At $\tau=m-2$, the same length already follows
from that bound for every $K\ge \tau+1$. The case $\tau=1$ recovers the
existence part of Theorem~\ref{thm:minimal-optimal}.

\begin{corollary}\label{cor:higher-minimal-five}
Let $\mathcal S_2=\{21,24,27,28,32,33,35,36,39,40,44\}$. For every
prime power $q$ and every $K$ in the additive semigroup generated by
$\mathcal S_2$, in particular for every $K\ge59$, there is a
nondegenerate $2$-minimal $[3K,K]_{q^5/q}$ code. It satisfies
$d_3\ge11$.
\end{corollary}
\begin{proof}
Here $(n,h,k,a,j)=(5,7,9,10,3)$, $\Delta=5$ and $T=3\ell/2$.
The stable profile imposes the six additional constraints
\[
 (1,1),\ (2,2),\ (3,4),\ (4,5),\ (5,7),\ (6,8).
\]
In Theorem~\ref{thm:simultaneous}, the dense pairs, with their
multiplicities, and the full-span terms of the additional constraints
have the following exponents in \eqref{eq:refined-counts} and
\eqref{eq:secondary-generic}:
\[
 \begin{array}{c|c|c|c}
 (d,c)&\text{mult.}&E_\ell(d,c)&E_6(d,c)\\ \hline
 (1,2)&2&1-2\ell&-11\\
 (1,3)&3&1-\ell/2&-2\\
 (2,5)&1&5-5\ell/2&-10\\
 (2,6)&2&4-\ell&-2
 \end{array}
 \qquad
 \begin{array}{c|c}
 (h_i,k_i)&\text{exponent of }q\text{ in \eqref{eq:secondary-generic}}\\ \hline
 (1,1)&2-s-2s\ell\\
 (2,2)&2-s-s\ell/2\\
 (3,4)&6-s-5s\ell/2\\
 (4,5)&5-s-s\ell\\
 (5,7)&12-s-3s\ell\\
 (6,8)&10-s-3s\ell/2
 \end{array}
\]
For every $s\ge7$ coprime to $5$ and every even $\ell\ge6$, we have
$T\ge9$, so $T\ge c$ for every dense pair and $sT\ge63\ge a_i$ for all
$i$; also $s\ell\ge h+1$. Bound the generic, dense and full-span terms
by \eqref{eq:refined-counts} and \eqref{eq:secondary-generic}, and use
\[
 R_{q,s,\ell}<q^{-5s\ell/2},\qquad
 \mathcal N_{q,s,\ell}<4q^{-s(\ell+1)}.
\]
The latter bound is from Proposition~\ref{prop:nondegenerate-kernel}.
All these bounds are nonincreasing in $q\ge2$, $s\ge7$ and even
$\ell\ge6$, by the displayed negative exponents and the monotonicity
of the product corrections. At $(q,s,\ell)=(2,7,6)$ their sum is
\[
 \begin{aligned}
 \mathfrak B_2={}&
 \frac{C_3(2^5)}{2(1-2^{-7})}
       \prod_{u=2}^{10}(1-2^{-u})^{-1}
 +\sum_{(d,c)}\mu_{d,c}\,2^{7E_6(d,c)}
                  C_c(2^7)C_d(2^{35})\\
 &+\sum_{i=1}^{6}\frac{2^{x_i}}{1-2^{-7}}
      \prod_{u=2}^{a_i}(1-2^{-u})^{-1}C_{j_i}(2^5)
   +2^{-105}+4\cdot2^{-49}<\frac{901}{1000},
 \end{aligned}
\]
where $\mu_{d,c}$ are the multiplicities in the first table and
$x_1,\ldots,x_6$ are the exponents in the second table at
$(s,\ell)=(7,6)$, namely $-89$, $-26$, $-106$, $-44$, $-121$ and $-60$. An explicit rational certificate
bounds the generic term by $900677/10^6$, the dense sum by $308/10^6$,
and the remaining terms together by $1/10^6$; the sum of these bounds is
$900986/10^6<901/1000$.
Consequently every such $(s,\ell)$ supplies a nondegenerate
kernel code of dimension $K=s\ell/2$, length $3K$, and with the
entire stable profile on its associated system.

Writing $\ell=2p$, the seed dimensions are the products $sp$ with
$s\ge7$, $\gcd(s,5)=1$ and $p\ge3$. The elements of $\mathcal S_2$ are
such products with $s\in\{7,8,9,11,13\}$ and $p\in\{3,4,5\}$. Their
semigroup has conductor $59$, as certified by the following consecutive
interval:
\[
 \begin{array}{c|c@{\qquad}c|c@{\qquad}c|c}
 K&\text{decomposition}&K&\text{decomposition}&K&\text{decomposition}\\ \hline
 59&24+35&66&27+39&73&33+40\\
 60&21+39&67&27+40&74&35+39\\
 61&21+40&68&24+44&75&35+40\\
 62&27+35&69&33+36&76&32+44\\
 63&24+39&70&35+35&77&33+44\\
 64&24+40&71&27+44&78&39+39\\
 65&21+44&72&28+44&79&35+44
 \end{array}
\]
Adding copies of the generator $21$ covers every $K\ge59$. The
preceding integer $58$ is neither a generator nor a sum of two
generators, and three generators total at least $63$; it is not a
seed dimension either, since $58=2\cdot29$. Below $59$, the semigroup
consists of
\[
 21,24,27,28,32,33,35,36,39,40,42,44,45,48,49,51,52,53,54,55,56,57.
\]
Direct sums preserve the profile by Lemma~\ref{lem:profile-sum}
and kernel nondegeneracy by Proposition~\ref{prop:nondegenerate-kernel}.
Proposition~\ref{prop:dictionary} gives $d_3\ge11$, and
\eqref{eq:higher-minimal-weight} proves $2$-minimality.
\end{proof}

\section{Maximum subspace designs and sum-rank codes}\label{sec:designs}

\subsection{The design condition}
Let $1\le h<r$. An \emph{$h$-design} in $\mathbb F_{q^n}^r$ is an ordered family
\[
 \mathcal U=(U_1,\ldots,U_t),\qquad U_i\le_{\mathbb F_q}\mathbb F_{q^n}^r,
\]
whose joint $\mathbb F_{q^n}$-span is $\mathbb F_{q^n}^r$ and such that
\begin{equation}\label{eq:design-definition}
 \sum_{i=1}^t\dimq(U_i\cap H)\le h
 \quad\text{for every }H\le_{\mathbb F_{q^n}} \mathbb F_{q^n}^r,\quad\dim_{\mathbb F_{q^n}} H=h.
\end{equation}
These are the finite-field subspace designs studied in \cite{SZ,SS}.
The broader notion, introduced in \cite{GX} and developed in
\cite{GK,GWX,GXY}, permits a separate bound on the sum of intersections
and has applications to list decoding.
The sum in \eqref{eq:design-definition} records the intersections
of the entire family; imposing a separate scatteredness condition
on every member does not suffice.

We will use the equivalent condition
\begin{equation}\label{eq:design-smaller}
 \sum_{i=1}^t\dimq(U_i\cap W)\le\dim_{\mathbb F_{q^n}} W
 \qquad(0\le\dim_{\mathbb F_{q^n}} W\le h).
\end{equation}
This follows from the usual extension argument. If the left side exceeds $d=\dim_{\mathbb F_{q^n}} W$, enlarge $W$ by a vector from some $U_i$ outside $W$; such a vector exists whenever $W\ne \mathbb F_{q^n}^r$ because the family spans. This increases the ambient dimension by one and increases the total intersection dimension by at least one. Repeating until dimension $h$ contradicts \eqref{eq:design-definition}. The converse is immediate.

Assume $n\ge h+1$. Every member of an $h$-design has dimension at most $rn/(h+1)$. Indeed, apply Theorem~\ref{thm:evasive-bound} inside its $\mathbb F_{q^n}$-span if that span has dimension greater than $h$; otherwise \eqref{eq:design-smaller} bounds its dimension by $h$, which is smaller than $rn/(h+1)$. A design is \emph{maximum} if
\[
 \dimq U_i=\frac{rn}{h+1}\qquad(1\le i\le t).
\]
Every member of a maximum design spans $\mathbb F_{q^n}^r$, by the same argument.

For $h=1$, the condition says that each $U_i$ is scattered and the linear sets
\[
 L_{U_i}=\{\langle u\rangle_{\mathbb F_{q^n}}:0\ne u\in U_i\}\subseteq\PG(r-1,q^n)
\]
are pairwise disjoint. This gives a concrete geometric interpretation of the construction below.

\subsection{An exact finite existence criterion}
Fix $1\le h<n$, $t\ge1$, and put
\begin{equation}\label{eq:Bt}
 B_{n,h,t}(q)=
 \sum_{\substack{a_1+\cdots+a_t=h+1\\0\le a_i\le n}}
       \prod_{i=1}^t\gauss{n}{a_i}{q}.
\end{equation}
The coefficient $\gauss{v}{u}{Q}$ counts $u$-dimensional subspaces of $\mathbb F_Q^v$.

For $t=1$, the scalar count below is the auxiliary-field
scattered-subspace argument of Giannoni and Neri~\cite{GN}.
The block formulation imposes the full design condition.

\begin{theorem}\label{thm:design-finite}
Let $\gcd(n,s)=1$, let $\ell$ be a positive integer with $(h+1)\mid n\ell$, and suppose $r=s\ell\ge h+1$. If
\begin{equation}\label{eq:design-finite}
 \frac{(q-1)B_{n,h,t}(q)}{(q^n-1)(q^s-1)}\le1,
\end{equation}
then there is a maximum $h$-design $(U_1,\ldots,U_t)$ in $\mathbb F_{q^n}^r$, with
\[
 \dimq U_i=N=\frac{rn}{h+1}.
\]
Under an $\mathbb F_{q^n}$-linear identification $\mathbb F_{q^n}^{\,r}\simeq\mathbb F_{q^{ns}}^{\,\ell}$, each $U_i$ can be chosen $\mathbb F_{q^s}$-linear.
\end{theorem}
\begin{proof}
Put $T=n\ell/(h+1)$, so that $N=sT$. Fix an $\mathbb F_q$-basis $\beta_1,\ldots,\beta_s$ of $\mathbb F_{q^s}$, also an $\mathbb F_{q^n}$-basis of $\mathbb F_{q^{ns}}$. Expand a vector $w\in \mathbb F_{q^{ns}}^T$ in this basis to obtain a vector in $\mathbb F_{q^n}^{sT}$, and denote its rank over $\mathbb F_q$ by $\rk_{q,\beta}(w)$.

Choose independent uniformly random matrices $A_i\in \mathbb F_{q^{ns}}^{\ell\times T}$ and put $A=[A_1|\cdots|A_t]$. Call a nonzero vector $w=(w_1,\ldots,w_t)\in \mathbb F_{q^{ns}}^{tT}$ \emph{bad} if
\[
 \sum_{i=1}^t\rk_{q,\beta}(w_i)\le h+1.
\]
We first show that avoiding bad vectors in $\ker_{\mathbb F_{q^{ns}}} A$ gives the desired family.

The restriction of each $A_i$ to $\mathbb F_{q^s}^T$ is injective: otherwise a nonzero kernel vector supported in that block would have $\rk_{q,\beta}=1$. Thus $U_i=A_i\mathbb F_{q^s}^T$ has $\mathbb F_q$-dimension $N$.

Suppose that an $\mathbb F_{q^n}$-space $H$ of dimension $d\le h$ violates \eqref{eq:design-smaller}. Select $d+1$ vectors from the intersections $U_i\cap H$, independent over $\mathbb F_q$ within each block. Write their inverse images as $z_{i,b}\in \mathbb F_{q^s}^T$. There is a nontrivial $\mathbb F_{q^n}$-linear relation among these $d+1$ vectors, with coefficients $c_{i,b}$. Define
\[
 w_i=\sum_b c_{i,b}z_{i,b}.
\]
The vectors $z_{i,b}$ in each block remain independent over $\mathbb F_{q^n}$, by the linear disjointness of $\mathbb F_{q^n}$ and $\mathbb F_{q^s}$. Hence some $w_i$ is nonzero. Its expanded coefficients lie in the $\mathbb F_q$-span of the corresponding $c_{i,b}$, so $\sum_i\rk_{q,\beta}(w_i)\le d+1\le h+1$. Since $Aw=0$, this contradicts the assumed avoidance.

To count bad vectors, choose coefficient spaces $D_i\le_{\mathbb F_q}\mathbb F_{q^n}$ of dimensions $a_i$ with $\sum_i a_i=h+1$. Vectors whose expanded coefficients in block $i$ lie in $D_i$ belong to $(\langle D_i\rangle_{\mathbb F_{q^s}})^T$. Linear disjointness gives $\dim_{\mathbb F_{q^s}}\langle D_i\rangle_{\mathbb F_{q^s}}=a_i$, so the total number of nonzero vectors covered by a fixed tuple is
\[
 q^{N(h+1)}-1=q^{ns\ell}-1.
\]
Every bad vector is covered by such a tuple, by enlarging its coefficient spaces if necessary. The total number of bad vectors is therefore at most
\[
 B_{n,h,t}(q)(q^{ns\ell}-1).
\]

The scalar group $\mathcal H_0=\mathbb F_{q^n}^*\mathbb F_{q^s}^*\le \mathbb F_{q^{ns}}^*$ has order
\[
 |\mathcal H_0|=\frac{(q^n-1)(q^s-1)}{q-1}.
\]
It preserves every $\rk_{q,\beta}(w_i)$ and acts freely on nonzero vectors. Indeed, multiplication by $\mathbb F_{q^n}^*$ scales coefficient spaces, while multiplication by $\mathbb F_{q^s}^*$ induces an invertible $\mathbb F_q$-linear change of the $\beta$-coordinates. Thus every $\mathbb F_{q^{ns}}$-line containing a bad vector contains at least $|\mathcal H_0|$ bad vectors. A fixed $\mathbb F_{q^{ns}}$-line lies in $\ker A$ with probability $q^{-ns\ell}$. The union bound gives
\[
 \Pr(\ker A\text{ contains a bad vector})
 \le \frac{B_{n,h,t}(q)}{|\mathcal H_0|}(1-q^{-ns\ell})<1
\]
under \eqref{eq:design-finite}.

It remains to check spanning. Let $r_i=\dim_{\mathbb F_{q^n}}\langle U_i\rangle_{\mathbb F_{q^n}}$. If $r_i\le h$, the already proved inequality \eqref{eq:design-smaller} gives $N\le r_i\le h$, whereas $N=rn/(h+1)\ge n>h$. If $r_i>h$, Theorem~\ref{thm:evasive-bound} inside this span gives $N\le r_i n/(h+1)$, hence $r_i\ge r$. Thus each $U_i$ spans $\mathbb F_{q^n}^r$, and the family is a maximum design.
\end{proof}

\begin{corollary}\label{cor:design-uniform}
A sufficient condition for \eqref{eq:design-finite}, simultaneously for all prime powers $q$, is
\begin{equation}\label{eq:design-uniform}
 s\ge h(n-1)+\log_2\left(2\binom{t+h}{h+1}4^{h+1}\right).
\end{equation}
\end{corollary}
\begin{proof}
The elementary Gaussian estimate $\gauss{v}{u}{q}<4q^{u(v-u)}$ and $a_i(n-a_i)\le a_i(n-1)$ imply
\[
 B_{n,h,t}(q)\le \binom{t+h}{h+1}4^{h+1}q^{(h+1)(n-1)}.
\]
There are at most $h+1$ nonzero parts in each composition; zero parts have Gaussian coefficient one. Also
\[
 |\mathcal H_0|=q^{n+s-1}\frac{(1-q^{-n})(1-q^{-s})}{1-q^{-1}}
       \ge \frac12q^{n+s-1}.
\]
It follows that $B_{n,h,t}(q)/|\mathcal H_0|\le1$ under \eqref{eq:design-uniform}, since $q\ge2$.
\end{proof}

\subsection{Every sufficiently large admissible ambient dimension}
\begin{lemma}\label{lem:design-sum}
If $(U_1,\ldots,U_t)$ and $(V_1,\ldots,V_t)$ are $h$-designs in $\mathbb F_{q^n}^{r_1}$ and $\mathbb F_{q^n}^{r_2}$, then
\[
 (U_1\oplus V_1,\ldots,U_t\oplus V_t)
\]
is an $h$-design in $\mathbb F_{q^n}^{r_1+r_2}$. If both designs are maximum, so is their direct sum.
\end{lemma}
\begin{proof}
For $W\le_{\mathbb F_{q^n}} \mathbb F_{q^n}^{r_1}\oplus \mathbb F_{q^n}^{r_2}$ with $\dim_{\mathbb F_{q^n}} W=d\le h$, let $\pi$ be projection to the first factor and put $W_0=W\cap \mathbb F_{q^n}^{r_2}$. Restricting this projection to each intersection gives
\[
 \dimq((U_i\oplus V_i)\cap W)
 \le \dimq(U_i\cap\pi W)+\dimq(V_i\cap W_0).
\]
Summing and using \eqref{eq:design-smaller} bounds the result by
$\dim_{\mathbb F_{q^n}}\pi W+\dim_{\mathbb F_{q^n}} W_0=d$. Spanning and the extremal dimensions are preserved by direct sums.
\end{proof}

\begin{theorem}\label{thm:design-all-r}
Fix integers $1\le h<n$ and $t\ge1$. There is an integer $R(n,h,t)$, independent of $q$, such that a maximum $h$-design with $t$ members exists in $\mathbb F_{q^n}^{\,r}$ for every prime power $q$ and every $r\ge R(n,h,t)$ satisfying $(h+1)\mid rn$.
\end{theorem}
\begin{proof}
Put $\nu=(h+1)/\gcd(n,h+1)$. Choose two distinct primes $s_1,s_2$ greater than $n$ and satisfying \eqref{eq:design-uniform}. Theorem~\ref{thm:design-finite}, with $\ell=\nu$, gives maximum designs in dimensions $\nu s_1$ and $\nu s_2$, for every $q$.

Every integer $R\ge(s_1-1)(s_2-1)$ can be written $R=us_1+vs_2$ with $u,v\ge0$. For example, choose $v\in\{0,\ldots,s_1-1\}$ with $vs_2\equiv R\pmod{s_1}$; then $R-vs_2\ge-(s_1-1)$ is a multiple of $s_1$, hence is nonnegative. Taking blockwise direct sums of the two designs proves the result for
\[
 r=\nu R,\qquad R\ge(s_1-1)(s_2-1).
\]
These are exactly all sufficiently large ranks allowed by the divisibility condition.
\end{proof}

\subsection{An exponential number of members}
The finite criterion also allows the number of members to grow with
the ambient dimension. The following estimate holds in every
sufficiently large admissible dimension.

\begin{lemma}\label{lem:balanced-coprime-parts}
Let $n\ge2$ and $M>6n$ be integers, and put
\[
 c_n(M)=
 \begin{cases}
 1,&\gcd(M,n)=1,\\
 2,&\gcd(M,n)>1\text{ and either $n$ is odd or $M$ is even},\\
 3,&\text{otherwise}.
 \end{cases}
\]
There are positive integers $s_1,\ldots,s_{c_n(M)}$, all coprime
to $n$, such that
\[
 M=\sum_i s_i,\qquad
 s_{\min}:=\min_i s_i>\frac{M}{c_n(M)}-n.
\]
If $n$ is odd or a power of two, then $c_n(M)\le2$ for every $M$.
\end{lemma}
\begin{proof}
If $\gcd(M,n)=1$, take $s_1=M$. To obtain two parts when $n$
is odd or $M$ is even, for every prime $p\mid n$ choose a residue
$a_p\not\equiv0,M\pmod p$. For $p>2$ at most two residues are
excluded. If $p=2$, the hypotheses imply that $M$ is even, so
only the residue zero is excluded. The Chinese remainder theorem
for these distinct primes supplies an integer $a$, which may be
chosen in $\{1,\ldots,n\}$, with
$\gcd(a,n)=\gcd(M-a,n)=1$. Put
\[
 s_1=a+n\left\lfloor\frac{M/2-a}{n}\right\rfloor,
 \qquad s_2=M-s_1.
\]
Both are coprime to $n$, and
$M/2-n<s_1\le M/2\le s_2<M/2+n$.

In the remaining case $n$ is even and $M$ is odd. Choose
\[
 s_1=1+n\left\lfloor\frac{M/3-1}{n}\right\rfloor.
\]
Then $s_1$ is coprime to $n$, is odd, and satisfies
$M/3-n<s_1\le M/3$. The remainder $M-s_1$ is even and at least
$2M/3$. Apply the preceding two-part construction to this remainder.
Each resulting part is greater than $(M-s_1)/2-n\ge M/3-n$.
All parts are positive because $M>6n$. Finally, if $n$ is a
power of two, every odd $M$ is coprime to $n$, so the third case
never occurs.
\end{proof}

\begin{theorem}\label{thm:design-growth}
Fix $1\le h<n$ and a prime power $q$, and put
$\nu=(h+1)/\gcd(n,h+1)$. Write $r=\nu M$ with $M>6n$, and
choose the parts in Lemma~\ref{lem:balanced-coprime-parts}.
Set
\begin{equation}\label{eq:design-number}
 t_0=\left\lfloor
       \frac{q^{(s_{\min}-h(n-1))/(h+1)}}{4\,2^{1/(h+1)}}
      \right\rfloor-h.
\end{equation}
If $t_0\ge1$, there is a maximum $h$-design with $t_0$ members
in $\mathbb F_{q^n}^{\,r}$. In particular, for every sufficiently
large $r$ with $(h+1)\mid rn$, there is such a design with
\begin{equation}\label{eq:design-exponential}
 t\ge q^{\,r/(c_n(M)\nu(h+1))-O_{q,n,h}(1)}.
\end{equation}
Uniformly in $M$, one may replace $c_n(M)$ by $2$ when $n$ is
odd or a power of two, and by $3$ otherwise. When
$\gcd(M,n)=1$, the coefficient $c_n(M)=1$ applies.
\end{theorem}
\begin{proof}
The estimates in the proof of Corollary~\ref{cor:design-uniform}
give the sufficient condition
\begin{equation}\label{eq:design-growth-criterion}
 2\,4^{h+1}\binom{t+h}{h+1}q^{h(n-1)-s}\le1
\end{equation}
for a maximum design in $\mathbb F_{q^n}^{\,\nu s}$, whenever
$\gcd(n,s)=1$ and $\nu s\ge h+1$.
Since $\binom{t+h}{h+1}\le(t+h)^{h+1}$, the value $t=t_0$
in \eqref{eq:design-number} satisfies this condition for every
$s\ge s_{\min}$. The lemma gives $s_{\min}>M/3-n>n$, so
$\nu s_i\ge h+1$ for every part. Thus the finite criterion,
with $\ell=\nu$, supplies a maximum design with $t_0$ members
in each dimension $\nu s_i$. Their direct sum, matching the
member indices, has ambient dimension $\nu\sum_i s_i=r$ and
is maximum by Lemma~\ref{lem:design-sum}.

To prove \eqref{eq:design-exponential}, let $X$ be the expression
inside the floor in
\eqref{eq:design-number}. For sufficiently large $M$ we have
$X\ge2(h+1)$, hence $t_0\ge X/2$. Since
$s_{\min}>M/c_n(M)-n$, this gives the explicit lower estimate
\[
 t_0\ge
 \frac{q^{(M/c_n(M)-n-h(n-1))/(h+1)}}{8\,2^{1/(h+1)}}.
\]
Taking base-$q$ logarithms proves \eqref{eq:design-exponential},
with a constant independent of $r$ and of the case defining
$c_n(M)$. The uniform alternatives follow from the last assertion
of the lemma.
\end{proof}

For every $1\le h<n$, condition \eqref{eq:design-smaller} on
one-dimensional $\mathbb F_{q^n}$-subspaces shows that the linear
sets of the members are scattered and pairwise disjoint. Writing
$N=rn/(h+1)$, each therefore has $(q^N-1)/(q-1)$ points in
$\PG(r-1,q^n)$, which has $(q^{nr}-1)/(q^n-1)$ points. This gives
\[
 t\le\left\lfloor
 \frac{(q^{nr}-1)(q-1)}{(q^n-1)(q^N-1)}
 \right\rfloor.
\]
For $h=1$, the theorem consequently gives exponentially many
pairwise disjoint maximum scattered linear sets in every sufficiently
large admissible ambient dimension. For $t\le q-1$, such families were
constructed in \cite[Theorem~4.11]{SS} whenever $rn$ is even, and the
first question of \cite[Problem~4.17]{SS} asks for more than $q-1$
members; Theorems~\ref{thm:design-all-r} and~\ref{thm:design-growth}
answer it for every sufficiently large admissible $r$. For $n=2$ and
odd $r$, the partitions into Baer subgeometries of
\cite[Theorem~4.7]{SS} have about $q^{r-1}$ members, many more than
Theorem~\ref{thm:design-growth} provides.

\subsection{The associated MSRD codes}
The sum-rank metric includes the rank and Hamming metrics and
admits both a matrix-space treatment \cite{BGR} and a support
theory for extension-field-linear codes \cite{MP19}.
Explicit MSRD families include linearized Reed--Solomon codes
\cite{MP18}, constructions from extended Moore matrices
\cite{MP22}, and the broader constructions of \cite{MP24}.
The geometry of such codes is developed in \cite{NSZ,SS}.
Linearized Reed--Solomon codes have at most $q-1$ blocks, each of
length at most $n$ (see \cite[Definition~2.9]{SS}); the codes obtained
below can have exponentially many blocks, each of length $N>n$.

For $x=(x^{(1)},\ldots,x^{(t)})\in(\mathbb F_{q^n}^N)^t$, define
\[
 \wt_{\rm SR}(x)=\sum_{i=1}^t\rk_q(x^{(i)}).
\]
A code in this block space is \emph{sum-rank nondegenerate} if the columns of a generator matrix are $\mathbb F_q$-independent within each block. This is a blockwise condition; it does not require independence of the entire collection of columns.

When $N\ge n$, an $\mathbb F_{q^n}$-linear code $\mathcal C\le(\mathbb F_{q^n}^N)^t$ of dimension $K$ and minimum sum-rank distance $d$ satisfies the matrix Singleton bound
\begin{equation}\label{eq:sr-singleton}
 nK\le N(tn-d+1).
\end{equation}
Indeed, expand into $t$ matrices of size $n\times N$ and delete a total of $d-1$ rows among the blocks. The projection is injective, since a nonzero difference in its kernel would have sum of ranks at most $d-1$. The remaining entries number $N(tn-d+1)$. A code attaining \eqref{eq:sr-singleton} is called MSRD, as in the long-block setting of \cite{BGR,SS}. In this setting duality preserves the MSRD property and the distances of nonzero dual codes sum to $tn+2$ \cite{BGR},
\cite[Theorem~2.2]{SS}.

\begin{corollary}\label{cor:msrd}
Let $(U_1,\ldots,U_t)$ be a maximum $h$-design in
$\mathbb F_{q^n}^{\,r}$, where $1\le h<n$ and $r\ge h+1$.
Assume $tn>h+1$, and write $N=rn/(h+1)$.
There exist sum-rank nondegenerate $\mathbb F_{q^n}$-linear MSRD
codes with $t$ equal blocks of length $N$ and parameters
\[
 \dim_{\mathbb F_{q^n}}\mathcal C=r,\qquad d_{\rm SR}(\mathcal C)=tn-h,
\]
and their duals satisfy
\[
 \dim_{\mathbb F_{q^n}}\mathcal C^\perp=tN-r,\qquad
 d_{\rm SR}(\mathcal C^\perp)=h+2.
\]
\end{corollary}
\begin{proof}
Choose an $\mathbb F_q$-basis of each design member as the columns of $G_i\in \mathbb F_{q^n}^{r\times N}$, and let $\mathcal C$ be the row space of $G=[G_1|\cdots|G_t]$. The spanning condition gives dimension $r$, and the columns in each block are $\mathbb F_q$-independent.

We claim $d_{\rm SR}(\ker G)\ge h+2$. Otherwise, choose a nonzero $c=(c^{(1)},\ldots,c^{(t)})\in\ker G$ of weight $b\le h+1$. In each block, a rank factorization writes
\[
 c^{(i)}=(\gamma_{i,1},\ldots,\gamma_{i,b_i})M_i,
\]
where $b_i=\rk_q(c^{(i)})$, the $\gamma_{i,a}$ are $\mathbb F_q$-independent, and the rows of $M_i\in\mathbb F_q^{b_i\times N}$ are independent. The $b_i$ vectors in $G_iM_i^{\mathsf T}$ are therefore independent over $\mathbb F_q$ within block $i$. All $b=\sum_i b_i$ such vectors satisfy the nontrivial $\mathbb F_{q^n}$-linear relation $Gc^{\mathsf T}=0$. Their $\mathbb F_{q^n}$-span has dimension at most $b-1\le h$, contradicting \eqref{eq:design-smaller}.

Applying \eqref{eq:sr-singleton} to the nonzero dual gives
\[
 n(tN-r)\le N(tn-d_{\rm SR}(\mathcal C^\perp)+1),
\]
so $d_{\rm SR}(\mathcal C^\perp)\le nr/N+1=h+2$. Thus the dual is MSRD. MSRD duality now gives $d_{\rm SR}(\mathcal C)=tn-h$, as claimed. This also recovers, in our notation, the design correspondence of
\cite[Theorem~3.16]{SS}.
\end{proof}

Applying the corollary to Theorem~\ref{thm:design-growth} gives,
for fixed $q,n,h$ and every sufficiently large admissible $r$,
MSRD codes whose number of equal long blocks satisfies
\eqref{eq:design-exponential}. The block length remains
$N=rn/(h+1)$, and the two distances are $tn-h$ and $h+2$.
Together with Theorem~\ref{thm:design-all-r}, the corollary thus gives
MSRD codes with block length $rn/(h+1)$ for every $1\le h<n$ and every
sufficiently large admissible $r$, as asked in
\cite[Problem~6.11]{SS}, although not by an explicit construction.

Astore and Zullo \cite{AZ} study length bounds for
$\mathbb F_{q^m}$-linear MSRD codes with each block length at
most $m$, attaining $d_{\rm SR}\le\text{length}-\dim+1$
\cite[Theorem~2.15 and Remark~2.16]{AZ}.
Here $m=n$ and $N=rn/(h+1)>n$ when $r>h+1$.
Our codes therefore lie in the long-block regime governed by
\eqref{eq:sr-singleton}, outside that hypothesis of their bounds.
In particular, the exponential block counts above do not address
the short-block extremal problems of \cite{AZ}.

\begin{remark}
The condition $tn>h+1$ excludes only the full-space endpoint $t=1$, $n=h+1$, where $\mathcal C=\mathbb F_{q^n}^r$, its distance is one, and its dual is zero. The maximum design theorem includes this endpoint; a minimum distance for the zero dual is not asserted.
\end{remark}

\begin{example}\label{ex:design-eight}
Take $q=2$, $n=4$, $h=1$, and $t=8$. Then
\[
 B_{4,1,8}(2)
 =8\gauss{4}{2}{2}+\binom82\gauss{4}{1}{2}^{\,2}
 =8\cdot35+28\cdot15^2=6580.
\]
With $s=9$ and $\ell=1$, the ratio in \eqref{eq:design-finite} is
\[
 \frac{6580}{(2^4-1)(2^9-1)}=\frac{188}{219}<1,
\]
whereas $t=9$ gives $561/511>1$; thus eight is the largest number of
members allowed by the criterion for these parameters.
Hence there are eight pairwise disjoint maximum scattered linear sets of
rank $18$ in $\PG(8,16)$. The associated code has eight blocks of length
$18$, $\mathbb F_{16}$-dimension $9$, and sum-rank distance $31$; its dual
has $\mathbb F_{16}$-dimension $135$ and distance $3$. These are existence
parameters certified by the counting argument, not an explicit list of
generating matrices.
\end{example}

\section{Further questions}\label{sec:discussion}

Efficient construction and verification remain natural questions
for the objects obtained here, especially when the number of
design members grows with the ambient dimension. The
scattered families of \cite{BGGM,BGM,BGGMN} and the explicit
design and MSRD constructions of \cite{GK,GXY,MP18,MP22,MP24,SS}
provide useful points of comparison.

The coprimality condition $\gcd(h,k+1)=1$ has two precise roles.
With $a=k+1$ and $j=a-h$, it makes the dense-stratum inequality
$jc\le ad$ strict for $0<d<j$, and it gives the boundary floor
identity used by the direct-sum profile. Without coprimality, our
bounds for strata with $jc=ad$ no longer decay as the matrix size
grows.
Controlling these strata would extend the evasive theorem and,
in particular, address ordinary minimal codes of length $2K$
over even extension degrees $m\ge4$, beyond the known short
cases \cite{BMN,XKH}. Proposition~\ref{prop:even-degree} makes this
precise for constructions by direct sums: all blocks but at most one
must be minimal MRD codes, and the MRD condition has the same critical
ratio $m/2$ as minimality. A single minimal MRD
$[2K',K']_{q^4/q}$ code with odd $K'\ge3$, for $q$ an odd power of $2$,
would already give $\varpi_{q^4/q}(K)=2K$ for all sufficiently large $K$
and those $q$: its direct sums with the codes of
Corollary~\ref{cor:quartic} are minimal MRD codes by
Proposition~\ref{prop:even-degree}(3), the integers $4$ and $K'$
generate every sufficiently large integer, and
Proposition~\ref{prop:minimal-lower} gives optimality.
This is a remaining part of the minimum-length problem in
\cite{ABNR,Scotti}. For higher minimality, the corresponding question
is whether $\gcd(m,\tau+1)=1$ can be removed.

For $\tau\ge2$, the construction proves
$\varpi_{q^m/q}(K,\tau)\le(\tau+1)K$ in the stated range.
A matching lower bound would require a new argument: the geometric
bound used for ordinary minimality does not apply to $k=m\tau-1$
when $\tau\ge2$. The case $m=5$, $\tau=2$ is a first concrete problem
within the higher-minimality framework of \cite{XKH}. There
$K+10\le\varpi_{q^5/q}(K,2)\le3K$ for every $K$ in the semigroup of
Corollary~\ref{cor:higher-minimal-five}, by \cite[Corollary~5.3(5)]{XKH}
and that corollary.

The stable profile determines the first two weights of the selected
ordinary minimal codes. Determining further generalized weights, or
prescribing simultaneous thresholds with equal critical ratios,
requires additional information. Idealizers and equivalence classes
would also help distinguish the resulting families beyond their
parameters; compare the weight analyses in \cite{MNT,ZZ,BGM}.
For designs, the gap between the exponential lower bound of
Theorem~\ref{thm:design-growth} and the point-packing bound remains
large in general; for $n=2$, $h=1$ and odd $r$ the packing bound is
attained by partitions into Baer subgeometries \cite[Theorem~4.7]{SS}.
Determining the largest number of members is the second question of
\cite[Problem~4.17]{SS}.

The auxiliary-field argument also invites comparison with the
existence questions for rank-metric codes over more general fields
studied by Neri and Zullo \cite{NZ}. The finite-field scalar count
used here does not directly settle those questions.
Within the finite-field setting, the threshold of
Corollary~\ref{cor:minimal-degree-five} is limited by the primary dense
pair $(1,2)$, whose exponent vanishes at $\ell=1$, and that of
Corollary~\ref{cor:higher-minimal-five} by the primary pairs $(1,3)$
and $(2,6)$, whose exponents vanish at $\ell=2$ and $\ell=4$. Smaller thresholds would require a finer count of dense
primary relations, or explicit seeds: by
Proposition~\ref{prop:profile-coding}, any minimal $[2K',K']_{q^5/q}$
code with $d_1=3$ can be added to the semigroup of
Corollary~\ref{cor:minimal-degree-five}, for the values of $q$ for which
it exists. For every $K\ge4$, Proposition~\ref{prop:minimal-lower} and
\cite[Corollary~5.4(1)]{XKH} give $2K\le\varpi_{q^5/q}(K)\le2K+2$, so the values of $K$
not covered by that corollary leave a gap of at most two.

\section*{Declaration of generative AI and AI-assisted technologies in the manuscript preparation process}

During the preparation of this work, the author used Claude (Anthropic)
in order to assist in drafting and revising parts of the text, including
the write-up of some proofs, and to check computations and references.
The author verified all arguments, numerical claims and citations. After
using this tool, the author reviewed and edited the content as needed and
takes full responsibility for the content of the published article.

\end{document}